\documentclass[a4paper,11pt]{article}
\usepackage[T1]{fontenc}
\usepackage[utf8]{inputenc}

\usepackage{lmodern}
\usepackage{amsmath}
\usepackage{amsthm}
\usepackage{amssymb}
\usepackage{authblk}
\usepackage{amsthm}
\usepackage{amssymb}
\usepackage{xparse}
\usepackage{thmtools}
\usepackage{stackrel}

\newcommand{\R}{\mathbb{R}}
\newcommand{\C}{\mathbb{C}}
\newcommand{\Z}{\mathbb{Z}}
\newcommand{\N}{\mathbb{N}}

\newcommand{\dd}{\mathrm{d}}
\newcommand{\CP}{\mathbb{CP}}
\renewcommand{\S}{\mathbb{S}}
\DeclareMathOperator{\Sp}{Sp}
\DeclareMathOperator{\tr}{tr}
\DeclareMathOperator{\dist}{dist}
\DeclareMathOperator{\diag}{diag}

\makeatletter
\newtheoremstyle{indented}
{7pt} 
{7pt} 
{} 
{1.5em} 
{\bfseries} 
{.} 
{.5em} 
{} 

\theoremstyle{definition}

\newtheorem{defn}{Definition}[section]

\theoremstyle{plain}
\newtheorem*{theorem*}{Theorem}

\newtheorem{theorem}{Theorem}

\newenvironment{preuve}{
	\noindent \textbf{Proof. }}{\hfill $\square$\medskip\par}

\newtheorem{prop}[defn]{Proposition}
\newtheorem{corr}[defn]{Corollary}

\newtheorem{lem}[defn]{Lemma}

\theoremstyle{definition}
\newtheorem{rem}[defn]{Remark} 

\title{Low-energy spectrum of Toeplitz operators: the case of wells}
\author{Alix Deleporte\thanks{deleporte@math.unistra.fr}}
\affil{Université de Strasbourg, CNRS, IRMA UMR 7501, F-67000
  Strasbourg, France}

\begin{document}

\maketitle
\let\thefootnote\relax\footnote{This work was supported by grant
  ANR-13-BS01-0007-01}

\begin{abstract}
In the 1980s, Helffer and Sj\"ostrand examined in a series of articles
the concentration of the ground state of a Schr\"odinger operator in
the semiclassical limit. In a similar spirit, and using the
asymptotics for the Szeg\H{o} kernel, we show a theorem about the
localization properties of the ground state of a Toeplitz operator,
when the minimal set of the symbol is a finite set of non-degenerate
critical points. Under the same condition on the symbol, for any
integer $K$ we describe the first $K$ eigenvalues of the operator.
\end{abstract}
\tableofcontents
\section{Introduction}
\subsection{Motivations}

In classical mechanics, the minimum of the energy, when it exists, is
a critical value, and any point in phase space achieving this minimum
corresponds to a stationary trajectory. In quantum mechanics, the
situation is quite different. A quantum state cannot be arbitrarily
localized in phase space and occupies at least some small amount of space, because of the uncertainty principle. Nevertheless, due to the correspondence principle, one expects the quantum states of minimal energy to concentrate, in some way, near the minimal set of the Hamiltonian, when the effective Planck constant is very small.

In a series of articles
\cite{helffer1984multiple,helffer1985multiple,helffer1985puits,helffer1986puits,helffer1987puits},
Helffer and Sj\"ostrand considered the Schr\"odinger operator
$P(\hbar)=-\hbar^2\Delta + V$, acting on $L^2(\R^n)$, where $V$ is a smooth
function. If $V$ is smooth, bounded from below and coercive, then the
infimum of the spectrum of $P(\hbar)$ is a simple eigenvalue. Helffer and
Sj\"ostrand then studied the concentration properties of
associated unit eigenvectors, named ground states, in the
semiclassical limit $\hbar\to 0$. It is well known that the ground state
is $O(\hbar^{\infty})$ outside any fixed neighbourhood of $\{x\in \R^n,
V(x)=\min(V)\}$. If there is only one such $x$, then the ground state
concentrates only on $x$. But if there are several minima, it is not
clear a priori whether the ground state is evenly distributed on them
or not.

In the first articles
\cite{helffer1984multiple,helffer1985multiple,helffer1985puits}, the
potential $V$ is supposed to reach its minimum only on a finite set of
non-degenerate critical points, named ``wells''. It is proven that
only some of the wells are selected by the ground state, that is, the
sequence of eigenvectors is $O(\hbar^{\infty})$ outside any fixed
neighbourhood of this subset of wells. The selected wells are the
flattest ones, in a sense that we will make clear later on. Sharper estimates lead to a control, outside of the wells, of the form $\exp(-C\hbar^{-1})$, where $C$ is expressed in terms of the Agmon distance to the selected wells. Moreover, when the potential $V$ has two symmetric wells, the ground state ``tunnels'' between these wells, so that there exists another eigenvalue which is exponentially close to the minimal one. 

In the two last papers \cite{helffer1986puits,helffer1987puits}, the potential $V$ is supposed to reach its minimal value on a submanifold of $\R^n$. Again, it is easy to prove that the ground state concentrates on this submanifold. From this fact, a formal calculus leads to the study of a Schr\"odinger operator, on the submanifold, with an effective potential that depends on the 2-jet behaviour of $V$ near the submanifold. The authors treated the case of an effective potential with one non-degenerate minimum, which they call the miniwell condition. In this case, the ground state concentrates only at the miniwell. On the contrary, when the minimal submanifold corresponds to a symmetry of $V$, the ground state is spread out on the submanifold.

This is an instance of what is called quantum selection: not all
points in phase space where the classical energy is minimal are
equivalent in quantum mechanics. When there is a finite set of minimal
points, only some of them are selected by the ground state. Similarly,
when the
minimal set is a submanifold, the ground state may select only one
point (or not). The series of articles
\cite{helffer1984multiple,helffer1985multiple,helffer1985puits}, and
also \cite{simon1983semiclassical,simon1984semiclassical} were
adapted to more general pseudodifferential operators
\cite{helffer1984puits,nakamura1998agmon,vu2001bohr}. In the physics literature, these effects are believed to appear in other settings; 
 for example, the miniwell condition was used in
 \cite{douccot1998semiclassical}, without mathematical justification,
 to study quantum selection effects for the Heisenberg model on spin systems, when the classical phase space is a product of 2-spheres. However, the arguments used by Helffer and Sj\"ostrand depend strongly on the fact that they deal with Schr\"odinger operators, when the phase space is $T^*\R^n$. Thus, it is a priori
not clear to which extent the quantum selection can be generalised to a quantization of compact phase spaces. 

We propose to study the K\"ahler quantization, which associates to a
symbol on a phase space a Toeplitz operator. In the particular case of
the coordinate functions on $\S^2=\CP^1$, the Toeplitz operators are
the spin operators \cite{borthwick2000introduction}, so that our
approach contains the physical case of spin systems. In this article
we prove (Theorem A) that quantum selection occurs in the case of wells, with
$O(\hbar^{\infty})$ remainder, extending some of the results of
\cite{helffer1984multiple,helffer1985multiple,helffer1985puits}. The
case of miniwells is in preparation. Exponential estimates will be the object of a separate investigation.

\subsection{K\"ahler quantization}
When a compact symplectic manifold is endowed with a K\"ahler structure, there is a natural way to define a quantization scheme, which is compatible with abstract geometric quantization \cite{kostant1970quantization,souriau1967quantification}. 

\begin{defn}
	A K\"ahler manifold is a complex manifold $M$ where the
        tangent space at each point is endowed with a hermitian
        metric, i.e. an inner product, whose
        imaginary part is a closed 2-form on $M$.
\end{defn}
In particular, a K\"ahler manifold has both a symplectic and a
Riemannian structure, which are respectively the imaginary part and
the real part of the inner product.

Let $L$ be a holomorphic complex line bundle over a K\"ahler manifold $M$, and $h$ denote
a hermitian metric on $L$. Let $\omega$ denote the imaginary part of
the hermitian metric. There exists a unique connection
(the Chern connection) compatible with $h$ and the complex structure. We wish to consider
a \emph{prequantum bundle} $(L,h)$, such that the curvature of the
Chern connection is $-i\omega$. This is
always locally possible, but the global existence of a prequantum bundle is
equivalent to the fact that $\omega/2\pi$ has integral cohomology
class (see for instance \cite{woodhouse1997geometric}, pp 158-162, or
Prop. 2.1.1 of \cite{kostant1970quantization}).
From now on we suppose that $\omega/2\pi$ has integral cohomology
class, and we let $(L,h)$ be a prequantum bundle.

There are two equivalent formulations for the basic objects of
K\"ahler quantization, one dealing with holomorphic sections of powers
of line bundles \cite{charles2000,charles2003berezin,le2014singular},
the other using equivariant functions on a circle bundle, the Grauert
tube \cite{BoutetdeMonvel1981,Zelditch2000,Shiffman2002}. In this
article we use the circle bundle approach.

Let $L^*$ denote the dual bundle of $L$, with $h^*$ the dual
metric. It is itself a Riemannian manifold. Define
\begin{align*}
D&=\{(m,v)\in L^*,\, h^*(v)<1\}\\X&=\partial D = \{(m,v)\in L^*,\,
h^*(v)=1\}.\end{align*}

Then $X$ is a circle bundle on $M$, with a $\mathbb{S}^1$-action on the fibres
$r_{\theta}:(m,v)\mapsto (m,e^{i\theta}v)$. We will also denote by
$\pi$ the projection from $X$ onto $M$.
As $X$ is a submanifold of a Riemannian manifold, it inherits a
volume form. Since $(L,h)$ is a prequantum bundle, $D$ is
pseudoconvex, and the volume form on $X$ coincides with the Levi form.
The scalar product on
$L^2(X)$ is related to the one on $L^2(M)$ via the $\S^1$
action. Indeed, if $s_0$ denotes any
smooth section of $X$, one has, for any $u,v\in L^2(X)$:
$$\langle u,v\rangle = \iint_{\mathbb{S}^1\times M} \overline{u}(r_{\theta}s_0(m))v(r_{\theta}s_0(m))
\dd \theta \dd m.$$

Let us consider the Hardy space, defined as follows. 
\begin{defn}
	\label{defn:Hardy-Szego}
	The \emph{Hardy space} on $X$, denoted by $H(X)$, is the closed subspace of $L^2(X)$ consisting of functions which are
	boundary values of holomorphic functions inside $D$.
	The orthogonal projector from $L^2(X)$ onto $H(X)$ is denoted by $S$,
	the \emph{Szeg\H{o}} projector.
\end{defn}

Using the $\S^1$ action, the space $H(X)$ can be further decomposed.
For $N\in
\N$, an element $f$ of $H$ is said to be $N$-equivariant when, for each $x\in X$ and
$\theta\in \S^1$, there holds $f(r_{\theta}x)=e^{iN\theta}f(x)$. The space of $N$-equivariant functions is denoted by $H_N(X)$; then $H(X)$
is the orthogonal sum of the different spaces $H_N(X)$ for $N\geq 0$. We will call
$S_N$ the orthogonal projection on $H_N(X)$. Then the Schwartz kernel
of $S_N$ is itself $N$-equivariant, that is: $$S_N(r_{\theta}x,r_{\phi}y)=e^{iN(\theta-\phi)}S_N(x,y).$$

For every $N$, the space $H_N(X)$ is finite-dimensional, the dimension
growing with $N$. To see this, note that the trace of $S_N$ is finite
as a consequence of Proposition \ref{prop:ma4.18}. It also comes from the fact
that $H_N(X)$ can also be formulated as a space of holomorphic
sections of an ample line bundle over the compact manifold $M$.

Now we define the \emph{K\"ahler quantization}, which associates to
any smooth function $f$ on $M$ a sequence of operators $(T_N(f))_{N\in
  \N}$:

\begin{defn}[Toeplitz operators]\label{def:toeplitz}
	Recall $\pi:X\to M$ is the natural projection. If $f\in C^{\infty}(M)$ is a smooth function, one defines the \emph{Toeplitz operator}
	with symbol $f$ as the sequence of operators $T_N(f):u \mapsto
	S_N(\pi^{*}\hspace{-0.2em}f\,u)$ from
	$H_N(X)$ to itself.
\end{defn}
In this article we are interested in the asymptotics, as $N\to
+\infty$, of Toeplitz operators and their eigenvectors.
Alternative conventions exist for the quantization (associating an operator to a symbol), though they define the same class of operators. The convention of Definition \ref{def:toeplitz} is sometimes called contravariant \cite{berezin1975general,charles2003berezin}. The reason for this choice is that we rely crucially on the positivity condition: if $f$ is real
and nonnegative, then $T_N(f)$ is nonnegative.

For any $N$, the operator $T_N(f)$
acts on a finite-dimensional space, moreover for real-valued $f$ this
operator is obviously
self-adjoint. Thus, the spectrum of $T_N(f)$ consists only of a finite
number of eigenvalues, each of which having a finite multiplicity. We will call the ``lowest eigenvalue'' the minimum of the spectrum of a Toeplitz operator.

We slightly extend the definition of Toeplitz operators in order to
deal with the K\"ahler manifold $M=\C^n$, which is not compact. This
does not affect the definitions of $H_N(M,L)$ and
$S_N$, except that the space $H_N(\C^n,L)$ has
infinite dimension in this case. If $f\in C^{\infty}(\C^n)$,
one can define the Toeplitz operator $T_N^{flat}(f)$ as an unbounded
operator, and it is an essentially self-adjoint operator with compact
resolvent, at least when the symbol is a positive quadratic form (see section 3.3).

\subsection{Main results}
In this article, we adapt the results from \cite{helffer1984multiple} to the setting of K\"ahler quantization. In particular, we are only interested in the following situation:
\begin{defn}\label{def:wellscond}
	A function $h\in C^{\infty}(M)$ is said to \emph{satisfy the
          wells condition} when the following is true:
\begin{itemize}
		\item $\min(h)=0$;
		\item Every critical point at which $h$ vanishes is non-degenerate.
	\end{itemize} 
\end{defn}

Observe that, by definition, Morse functions whose minimum is zero satisfy the wells condition, as does
the square modulus of a generic holomorphic section of $L^{\otimes N}$ for $N$ large. Note that a function
that satisfies the wells condition has a finite cancellation set.

We need the following definition to state our main theorems:
\begin{defn}\label{defn:precon}
	Let $Z$ be a subset of $M$, and let
	$$V_{\delta}(N)=\{(m,v)\in X,\,\dist(m,Z)>N^{-\delta}\}.$$
	
	A sequence $(u_N)_{N\in \N}$ of norm 1 functions in $L^2(M,L)$  is said to \emph{concentrate} on $Z$ when, for every $\delta \in [0,\frac 12)$,  one has $$\|u_N1_{V_{\delta}(N)}\|_{L^2(X)}=O(N^{-\infty}).$$
\end{defn}

Note that concentration, in the sense of the definition above, implies
microsupporting in the sense of Charles \cite{charles2003berezin},
that is, for any open set $V$ at positive distance from $Z$, as $N\to +\infty$,
one has
$\|u_N1_{V}\|_{L^2}=O(N^{-\infty})$. The microsupport is contained in
the concentration set, while the concentration set is included in any
open neighbourhood of the microsupport.

In Subsection 2.2, we consider convenient local maps of ``normal coordinates'' around any point $P\in M$, which preserve infinitesimally the K\"ahler structure.
If a non-negative function $h$ vanishes with positive Hessian at $P\in M$, the 2-jet of $h$ at $P$ reads in these coordinates as
a positive quadratic form $q(P)$ on $\C^n$. The first
eigenvalue $\mu$ of the Toeplitz operator $T_1^{flat}(q(P))$ (which we
call \emph{model quadratic operator}) does not depend on the choice of
normal coordinates. We define this value to be $\mu(P)$.

Let now $h$ be a smooth function on $M$ that satisfies the wells condition of Definition \ref{def:wellscond}. 
\begin{theorem}
	For every $N\in \N$, let $\lambda_N$ be the
	first eigenvalue of the operator $T_N(h)$, and $u_N$ an associated
	normalized eigenfunction.
	Then the sequence $(u_N)_{N\in \N}$ concentrates on the vanishing points of $h$ on
	which $\mu$ is minimal. 
	
	If there is only one such point $P_0$, then
	there is a real sequence $(a_k)_{k\geq 0}$ with $a_0=\mu(P_0)$ such that, for each $K$, one has
	$$\lambda_N=N^{-1}\sum_{k=0}^KN^{-k}a_k+O(N^{-K-2}).$$ Moreover, $\lambda_N$ is simple, and there
	exists $C>0$ such that
	$\lambda_N$ is the only eigenvalue of $T_N$ in the interval $[0,N^{-1}(\mu(P_0)+C)]$.
\end{theorem}

\begin{theorem}
 Let $C>0$. There is a bounded number of eigenvalues (counted with multiplicity) of $T_N(h)$ in the interval $[0,CN^{-1}]$. 
	More precisely, for $C'>C$, let $K$ and $(b_k)_{1\leq k \leq K}$ be such
        that $$\{b_k,\,k\leq K\}=\bigcup_{\stackrel{P\in M}{h(P)=0}}\Sp\left(T_1^{flat}(q(P))\right) \cap [0,C']$$ with multiplicity.
	Then one can find $c>0$ and a list of real numbers $(c_k)_{1\leq k \leq K}$ such that, for each $k$, one of the eigenvalues of $T_N(h)$ lies in the interval $$[N^{-1}b_k+ N^{-3/2}c_k - cN^{-2},N^{-1}b_k+ N^{-3/2}c_k + cN^{-2}].$$
	Moreover, there are at most $K$ eigenvalues of $T_N(h)$ in
        $[0,CN^{-1}]$ and each of them belongs to one of the intervals
        above.

        Among the smooth functions satisfying the wells condition,
        there is a dense open subset of ``non-resonant'' symbols such
        that, for every $k\geq 0$, the $k$-th eigenvalue of the
        associated operator has an
        asymptotic expansion in powers of $N^{-1/2}$.
\end{theorem}

The case of ``miniwells'', a transposition of \cite{helffer1986puits},
will be treated in future work. Under analyticity conditions, we also
hope to state results on exponential decay in the forbidden region, as
in \cite{helffer1985puits}.
In the one-dimensional case, a full asymptotic expansion for the first eigenvalues of
$T_N(h)$ was given in \cite{le2014singular}, with a fixed domain of
validity $[0,E_0]$. 

\subsection{Methods -- semiclassical properties of K\"ahler quantization}

If $M$ is compact or $\C^n$, the K\"ahler quantization has many similarities with the Weyl
quantization on cotangent bundles. One can indeed find a star product on the space of
formal series $C^{\infty}(M)[[\eta]]$ that coincides with the
composition of Toeplitz operators when $\eta=N^{-1}$ (see
\cite{schlichenmaier2000deformation} for a short proof), that is, such
that $$T_N(f\star g)=T_N(f)T_N(g) + O(N^{-\infty}).$$

Thus, the limit $N\to +\infty$ for Toeplitz operators can be thought
of as a semiclassical limit, with semiclassical parameter $N^{-1}\to 0$. Unless otherwise stated, we will state
results under this limit.

It has been known since at least \cite{charles2003berezin} that
there is a microlocal equivalence between the semiclassical calculus of
Toeplitz
operators and that of Weyl quantization. Such a correspondence was
already given in the homogeneous setting (without a semiclassical parameter) in \cite{BoutetdeMonvel1981}.

Thus, a possible approach to the spectral study of Toeplitz operators
(such as this one, which focuses on low-lying eigenvalues) would be
a conjugation by a Fourier Integral Operator to an operator known by
previous work. This could be a pseudodifferential operator, in the
spirit of \cite{helffer1984multiple}, or a Toeplitz operator with a
simpler symbol, cf
\cite{vu2000formes,raymond2013geometry}. However, each of these
approaches require a priori results on the concentration of eigenvectors. 

We will use a direct approach in this paper. Indeed, our future work
(in preparation) will focus on
the case when the minimal set of the symbol is a submanifold, where a
priori concentration is not known, so it is unclear whether the previous
approaches are sufficient. Moreover, the main theorems in this paper
depend on subprincipal effects, and the criterion for quantum
selection would be less natural if we should keep track of it through
a Fourier Integral Operator. Finally, we believe that Proposition
\ref{prop:loc} is of independent interest. It can easily be generalized into a result on the
microsupport of low-energy states for any smooth symbol, and it does not depend on
estimates on the asymptotics of the Szeg\H{o} kernel but only on the
nature and symbolic calculus of Toeplitz operators. It could be used as an elementary proof of microsupporting
for pseudodifferential operators.

\subsection{Outline}
We review in Section 2 the definitions and semiclassical properties of
the Szeg\H{o} kernel. Using well-known results about its semiclassical
expansion \cite{Shiffman2002,charles2003berezin,ma2006, Berman2008},
we derive Proposition \ref{prop:universal}, which states that the
Szeg\H{o} kernel on $\C^n$ is a local model for any Szeg\H{o} kernel.

In Section 3, we remind the reader of the symbolic properties of Toeplitz operators \cite{schlichenmaier2000deformation}. The state of the art is such that one can compose Toeplitz operators 
with classical symbols. We then show, with a new method, a standard result on localization: low-energy eigenvalues concentrate where the symbol is minimal. Finally, we study in detail a particular case of Toeplitz operators, when the base manifold is $\C^n$ and the symbol is a positive quadratic form.

Section 4 is devoted to the proof of Theorem A. We build an
approximate eigenfunction of the Toeplitz operator and prove that the
corresponding eigenvalue is the lowest one. The most important part is
Proposition \ref{prop:well-ex} for which we use the same method than
in \cite{helffer1984multiple}. For this, we consider the Hessian of
$h$ at a cancellation point, as read in local coordinates; this is a
real quadratic form $q$ on $\C^n$. Then we compare the Toeplitz operator $T_N(h)$ with the Toeplitz operator $T_N^{flat}(q)$, which we call \emph{model quadratic operator}.

In Section 5, we modify the arguments used in Section 4 to describe,
under the same hypotheses on the symbol, the spectrum of a Toeplitz
operator in the interval $[0,CN^{-1}]$ where $C>0$ is arbitrary
(Theorem B).

The Appendix is independent from the two main results of the paper. We
recover, in the K\"ahler setting, the off-diagonal estimate for the Szeg\H{o} kernel of
\cite{ma2006, charles2003berezin, Berman2008}, in local
coordinates. For this we use the techniques developed in
\cite{Zelditch2000, Shiffman2002}, which yield estimates on a
shrinking scale, and slightly modify them to recover an estimate on a
fixed scale.

\setcounter{theorem}{0}

\section{The Szeg\H{o} projector}

\subsection{Bargmann spaces}
As a helpful illustration for the general case (which originates from
\cite{bargmann1961hilbert,bargmann1967hilbert}, see also
\cite{Folland1989harmonic}, pp. 39-51), we first consider the usual $n$-dimensional
complex space $\C^n$, with the natural K\"ahler
structure, with $\omega=\sum_{i=1}^n \dd z_i\wedge \dd \overline{z}_i$. In this example, the Szeg\H{o} kernel is explicit.

Because $\C^n$ is contractible, the bundle $L$ is isomorphic to
$\C^{n+1}$, but the hermitian structure $h$ is not the flat one,
for which the associated curvature is zero. Indeed, one can show that
$h(m,v)=e^{-|m|^2}|v|^2$ is the correct choice. Here, the spaces $H_N(\C^n,L)$ are called the
\emph{Bargmann spaces} and will be denoted by $\mathcal{B}_N$. They can be expressed as
$$\mathcal{B}_N = L^2(\C^n) \cap \left\{z \mapsto e^{-\frac N2
  |z|^2}f(z),\,f\text{ holomorphic in $\C^n$ }\right\}.$$

The space $\mathcal{B}_N$ is a closed subspace of the Hilbert space
$L^2(\C^n)$ and inherits its scalar product: $$\langle f,g\rangle = \int_{\C^n}\overline{f}g.$$

The functions in $\mathcal{B}_N$ are not holomorphic for the standard
structure. However, let us introduce the following deformation of
$\overline{\partial}$:
$$\overline{\mathfrak{d}}_N=e^{-\frac N2 |z|^2}\overline{\partial} e^{\frac N2 |z|^2}
= \overline{\partial} + \frac N2 z.$$
We will further denote by $\overline{\mathfrak{d}}_{Ni}$ the $i$-th component
of $\overline{\mathfrak{d}}_N$. The space $\mathcal{B}_N$ is the space of $L^2$ functions in the kernel of $\overline{\mathfrak{d}}_N$. The adjoint of $\overline{\mathfrak{d}}_N$ is $\mathfrak{d}_N=e^{-\frac N2 |z|^2}\partial e^{\frac N2 |z|^2}$.
The orthogonal projector on $\mathcal{B}_N$ has a Schwartz
kernel. Indeed, one Hilbert basis of $\mathcal{B}_N$ is the family
$(e_{\nu})_{\nu \in \Z^n}$ with
$$e_{\nu}(z)=N^n\cfrac{N^{|\nu|/2}z^{\nu}}{\pi^n\sqrt{\nu!}}e^{-\frac
  N2 |z|^2}.$$ Hence the kernel may be expressed as:
\begin{equation}\label{eq:projbarg}\Pi_{N}(x,y)=\cfrac{N^n}{\pi^n}\exp\left(-\frac
  N2|x|^2-\frac N2|y|^2+Nx\cdot
\overline{y}\right).\end{equation}
Note that, by definition, $\Pi_N$ commutes with
$\overline{\mathfrak{d}}_N$. Moreover $$[\Pi_N, \overline{z}_i]= \Pi_N \mathfrak{d}_{Ni}.$$

The space $\mathcal{B}_N$ is isometric to $\mathcal{B}_1$ by an isometric
dilatation (or scaling) of factor $N^{1/2}$:
\begin{align*}
  \mathcal{B}_N & \leftrightarrow \mathcal{B}_1\\
f &\mapsto N^{-n}f(N^{-1/2}\cdot).
\end{align*}

Moreover, there is a unitary transformation between $\mathcal{B}_1$ and
$L^2(\R^n)$, called the \emph{Bargmann transform}. The transformation
$B_1:L^2(\R^n)\mapsto \mathcal{B}_1$ reads:
$$B_1f(z)=e^{-\frac 12 |z|^2}\int \exp\left[-\left(\frac 12 z\cdot z+ \frac 12 x \cdot x-\sqrt{2}z \cdot
x\right)\right]f(x)\dd x.$$ This transformation conjugates the position operators
$z_i$ into the position operators $x_i$, and the momentum operators
$\mathfrak{d}_{1\,i}$ into the momentum operators $\frac 1i\frac{\partial}{\partial x_i}$.

From $B_1$, one can deduce an isometry from $\mathcal{B}_N$ to $L^2(\R^n)$ by composing the scaling isometry and the Bargmann transform.

One noteworthy subspace of $\mathcal{B}_1$ is the dense subset of
functions $f\in \mathcal{B}_1$ such that $fP\in \mathcal{B}_1$ for any
polynomial $P\in \C[z]$. This space is denoted by $\mathcal{D}$. Any element of
the previously given Hilbert basis belongs to $\mathcal{D}$ and the
Bargmann transform is a bijection from $\mathcal{S}(\R^n)$ to
$\mathcal{D}$; the preimage of $e_{\nu}$ is the function $x\mapsto
C_{\nu}x^{\nu}e^{-|x|^2/2}$, where $C_\nu$ is a normalizing factor. 
Moreover, because of the commutation relations above, the image of
$\mathcal{S}(\C^n)$ by the Szeg\H{o} projector $\Pi_N$ is $\mathcal{D}$.
\subsection{Semiclassical asymptotics}
Semiclassical expansions of $S_N$ are
derived in
\cite{Zelditch2000,Shiffman2002, ma2007holomorphic, charles2000, Berman2008}, in different
settings. In \cite{Zelditch2000,Shiffman2002}, the Fourier Integral Operator approach is used to prove an asymptotic expansion of $S_N$ in a neighbourhood of size $N^{-1/2}$ of a point. In \cite{charles2000,ma2007holomorphic,Berman2008}, one derives asymptotic expansions of $S_N$ in a neighbourhood of fixed size of a point.

The Szeg\H{o} kernel is rapidly decreasing away from the diagonal as $N\to +\infty$:

\begin{prop}[\cite{charles2003berezin}, Corollary 1, or \cite{ma2006}, Prop. 4.1 in a more general setting]\label{prop:ma4.1}
	For every $k\in \N$ and $\epsilon>0$, there exists $C>0$ such that, for every $N\in \N$, for every $x,y\in X$, if $$\dist(\pi(x),\pi(y))\geq \epsilon,$$ then
	$$|S_N(x,y)|\leq CN^{-k}.$$
\end{prop}

The analysis of the Szeg\H{o} kernel near the diagonal requires a convenient choice of coordinates.
Let $P_0\in M$. The real tangent space $T_{P_0}M$ carries a Euclidian structure and an almost complex structure coming from the K\"ahler structure on $M$. We then can (non-uniquely) identify $\C^n$ with $T_{P_0}M$.

\begin{defn}\label{defn:normap}
	Let $U$ be a neighbourhood of $0$ in $\C^n$ and $V$
	be a neighbourhood of $P_0$ in $M$. Let $\pi$ denote the
        projection from $X$ to $M$.
        Let $\R$ cover $\S^1$. The group action
        $r_{\theta}:\S^1\to X$ lifts to a periodic action from
        $\R$ to $X$, which we will also call $r_{\theta}$.
	A smooth diffeomorphism $\rho:U\times \R \to \pi^{-1}(V)$ is said to be a
	\emph{normal map} or map of \emph{normal coordinates} under
	the following conditions:
	\begin{itemize}
		\item $\forall z \in U,\,\forall \theta \in \R,\,\rho(z,\theta)=r_{\theta}\rho(z,0)$;
		\item Identifying $\C^n$ with $T_{P_0}M$ as previously, one has: $$\forall (z,\theta)\in U\times \R,\,\pi(\rho(z,\theta))=\exp(z).$$ 
	\end{itemize}
\end{defn}
Through this paper we will often read the kernel of $S_N$ in normal
coordinates. Let $P_0\in X$ and $\rho$ a normal map on $X$ such that
        $\rho(0,0)=P_0$. For $z,w\in \C^n$ small enough and $N\in \N$,
        let $$S^{P_0}_N(z,w):=
        e^{-iN(\theta-\phi)}S_N(\rho(z,\theta),\rho(w,\phi)),$$ which
        does not depend on $\theta$ and $\phi$ as $S_N$ is $N$-equivariant.
	
	The following proposition states that, as $N\to +\infty$, in normal coordinates, the Szeg\H{o} kernel has an asymptotic expansion whose first term is the flat kernel of equation (\ref{eq:projbarg}):

\begin{prop}[\cite{ma2006}, theorem 4.18]\label{prop:ma4.18}
	There exist $C>0$, $C'>0$, $m\in \N$, $\epsilon>0$ and a
        sequence of polynomials $(b_j)_{j\geq 1}$, with $b_j$ of same
        parity as $j$, such that, for any $N\in \N$, $K\geq 0$ and $|z|,|w|\leq \epsilon$, one has:
	\begin{multline}\label{eq:SNexpfixed}
		\left|S^{P_0}_N(z,w)-\Pi_N(z,w)\left(1+\sum_{j=1}^KN^{-j/2}b_j(\sqrt{N}z,\sqrt{N}w)\right)\right|\leq \\CN^{n-(K+1)/2}\left(1+|\sqrt{N}z|+|\sqrt{N}w|\right)^me^{-C'\sqrt{N}|z-w|}+O(N^{-\infty}).
	\end{multline}
\end{prop}

\begin{rem}
We will use Proposition \ref{prop:ma4.18} as a black box, as we do
not want to divert the reader into considerations on the asymptotics
of $S_N$, which are more technical than the rest of this paper.

The scope of \cite{ma2006} is much more general than the case of
K\"ahler manifolds; by specialising to this case, one obtains stronger
estimates. Indeed, a result very close to this proposition can be found in \cite{Berman2008}, and also in \cite{charles2003berezin}, Theorem
2. However, these
results are stated without local coordinates, hence the link with the
Bargmann spaces is not obvious. 

For the sake of the argument, we derive in the Appendix a precised
formulation of this stronger version, adapting the techniques
presented in \cite{Shiffman2002}.
\end{rem}

\begin{rem}
	The Proposition \ref{prop:ma4.18} gives asymptotics for the
        kernel of $S_N$, read in local coordinates. However, the
        normal maps of Definition \ref{defn:normap} do not preserve
        the volume form, except infinitesimally on the fibre over
        $P_0$. For the associated operators to be preserved, one has to
        pull-back Schwartz kernels as half-forms. We claim that it
        does not change the structure of the asymptotics.
	
	Indeed, if $\dd \text{Vol}$ is the volume form on $X$ and $\dd Leb$
        is the Lebesgue form on $\C^n$, one has, for any normal map
        $\rho$: $$\rho^*(\dd Leb \otimes \dd \theta)= a \,\dd \text{Vol},$$
        for some function $a$ on the domain of $\rho$ with $a(0)=1$. We want to study
        the asymptotics of $(z,w)\mapsto
        S_N^{P_0}(z,w)\sqrt{a(z)a(w)}$, which is the kernel of the
        pull-back of $S_N$.

 The function $(z,w)\mapsto
        \sqrt{a(z)a(w)}$ is smooth on the domain of $\rho$. We write
        the Taylor expansion of this function at $0$ as:
$$\sqrt{a(z)a(w)}=1+\sum_{j=1}^Ka_j(z,w)+O(|z|^{K+1},|w|^{K+1})$$ where $a_j$ is
homogeneous of degree $j$, so that
$a_j(z,w)=N^{-j/2}a_j(\sqrt{N}z,\sqrt{N}w)$.

We let now $\widetilde{b}_j$ be such
that \begin{multline*}\left(1+\sum_{j=1}^KN^{-j/2}b_j(\sqrt{N}z,\sqrt{N}w)\right)\left(1+\sum_{j=1}^KN^{-j/2}a_j(\sqrt{N}z,\sqrt{N}w)\right)\\=1+\sum_{j=1}^KN^{-j/2}\widetilde{b}_j(\sqrt{N}z,\sqrt{N}w)+O(N^{-(K+1)/2}).\end{multline*}

Then
\begin{multline*}
		\left|S^{P_0}_N(z,w)\sqrt{a(z)a(w)}-\Pi_N(z,w)\left(1+\sum_{j=1}^KN^{-j/2}\widetilde{b}_j(\sqrt{N}z,\sqrt{N}w)\right)\right|\leq \\CN^{n-(K+1)/2}\left(1+|\sqrt{N}z|+|\sqrt{N}w|\right)^me^{-C'\sqrt{N}|z-w|}+O(N^{-\infty}).
	\end{multline*}

	Hence, the effects of the volume form can be absorbed in the error terms of equation (\ref{eq:SNexpfixed}), and the Proposition \ref{prop:ma4.18} also holds when $S_N$ is replaced by the corresponding half-form.
\end{rem}

Thus, we can use the asymptotics of Proposition \ref{prop:ma4.18} to
study how the operator $S_N$ acts. For instance, we are able to refine the Proposition \ref{prop:ma4.1}:

\begin{corr}\label{corr:localization}
		For every $k\in \N$ and $\delta\in [0,1/2)$, there exists $C>0$ such that, for every $N\in \N$, for every $x,y\in X$ with $\dist(\pi(x),\pi(y))\geq N^{-\delta}$, one has:
		$$|S_N(x,y)|\leq CN^{-k}.$$
		
		In particular, if $u\in L^2(X)$ is $O(N^{-\infty})$ outside the pull-back of a ball of size $N^{-\delta}$, then $S_N(u)$ is $O(N^{-\infty})$ outside the pull-back of a ball of size $2N^{-\delta}$.
\end{corr}

\subsection{Universality}
In the previously given local expansions of the Szeg\H{o} kernel (\ref{eq:SNexpfixed}), the
dominant term is the projector on the Bargmann spaces of equation (\ref{eq:projbarg}). Thus the
Bargmann spaces appear to be a universal model for Hardy spaces, at
least locally. To make this intuition more precise, we derive a useful
proposition.

We can pull-back by a normal map the kernel of the projector $\Pi_N$ by the following formula:
$$\rho^*\Pi_N(\rho(z,\theta),\rho(w,\phi)):=e^{iN(\theta-\phi)}\Pi_N(z,w).$$

By convention, $\rho^*\Pi_N$ is zero outside $\pi^{-1}(V)^2$.

\begin{prop}[Universality]\label{prop:universal}Let
  $\epsilon>0$. There exists $\delta\in (0,1/2)$, a constant $C>0$ and
  an integer $N_0$ such that, for any $N\geq N_0$, for any function
  $u\in L^2(X)$ whose support
is contained in the fibres over a ball on $M$ of radius $N^{-\delta}$, one has:
   $$\|(\rho^*\Pi_N)u-S_Nu\|_{L^2(X)}\leq CN^{-1/2+\epsilon}\|u\|_{L^2(X)}.$$
   
\end{prop}
\begin{preuve}
Let again $S_N^{P_0}: (z,\theta, w, \phi) \mapsto
e^{-iN(\theta-\phi)}S_N(\rho(z,\theta),\rho(w,\phi))$ denote the
kernel $S_N$ as read in local coordinates, which does not in fact
depend on $(\theta,\phi)$. 

  Equation (\ref{eq:SNexpfixed}), for $K=0$, can be formulated as:
  \begin{equation}\label{eq:SNexpfixed0}
    S_N^{P_0}(z,w)=\Pi_N(z,w)+R(z,w)+O(N^{-\infty}),
  \end{equation}
  with $$|R(z,w)|\leq
  CN^{n-1/2}(1+|\sqrt{N}z|+|\sqrt{N}w|)^me^{-C'\sqrt{N}|z-w|}$$ for
  every $z$ and $w$ such that $(z,0)$ and
$(w,0)$ belong to the domain of $\rho$.
  
Let $\delta \in (0,1/2)$ and $u$ a function contained in the pull-back of a ball of size $N^{-\delta}$.

Let $v=S_Nu-(\rho^*\Pi_N)u$. Because of Corollary 
\ref{corr:localization}, $v$ is $O(N^{-\infty})$ outside
$\rho(B(0,4N^{-\delta})\times \S^1)$. Hence, up to a $O(N^{-\infty})$
error, it is sufficient to control the kernel of $S_N-\rho^*\Pi_N$ on
$\rho(B(0,4N^{-\delta})\times \S^1)\times \rho(B(0,4N^{-\delta})\times
\S^1)$, where equation (\ref{eq:SNexpfixed0}) is valid.

It remains to estimate the norm of the operator with kernel $R$, using a standard result of operator theory:
\begin{lem}[Schur test]
	Let $k \in C^{\infty}(V\times V)$ be a smooth function of two variables in an open subset $V$ of $\R^d$. Let $K$ be the associated unbounded operator on $L^2(V)$.
	
	Let $$\|k\|_{L^{\infty}L^1} := \max\left(\sup_{x\in V}\|k(x, \cdot)\|_{L^1(V)},\sup_{y\in V}\|k(\cdot, y)\|_{L^1(V)}\right).$$ If $\|k\|_{L^{\infty}L^1}$ is finite, then $K$ is a bounded operator. Moreover $$\|K\|_{L^2(V)\mapsto L^2(V)}\leq \|k\|_{L^{\infty}L^1}.$$
\end{lem}

Thus, we want to estimate the quantity:
$$\sup_{|z|\leq 4N^{-\delta}}\int_{|w|\leq 4N^{-\delta}}N^{n-1/2}(1+|\sqrt{N}z|+|\sqrt{N}w|)^me^{-C'|z-w|}.$$
After a change of variables and up to a multiplicative constant, it remains to estimate:
$$N^{-1/2}\sup_{|z|\leq 4N^{1/2-\delta}}\int_{|u|\leq 4N^{1/2-\delta}}\left(1+|z|+|u|\right)^me^{-C|u|}.$$
This quantity is $O(N^{(m-1)\frac 12 - m \delta})$. Thus, for any $\epsilon>0$, there exists $\delta$ such that the above quantity is $O(N^{-\frac 12+\epsilon})$.

By the Schur test, the $L^2$ norm of a symmetric kernel operator is
controlled by the $L^{\infty}L^1$ norm of the kernel. When restricted
on $B(0,4N^{-\delta})^2$, the kernel of $S_N^{P_0}-\Pi_N$ has a
$L^{\infty}L^1$ norm of order $N^{-\frac 12 + \epsilon}$, from which we can conclude.
\end{preuve}

\section{Toeplitz operators}
\subsection{Calculus of Toeplitz operators}
%
%

The composition of two Toeplitz operators is a formal series of Toeplitz operators. The theorem 2.2 of
\cite{schlichenmaier2000deformation} states for instance that there
exists a formal star-product on $C^{\infty}(M)[[\eta]]$, written as $f\star g =
\sum_{j=0}^{+\infty}\eta^jC_j(f,g)$, that coincides with the Toeplitz
operator composition: as $N\to +\infty$, one has, for every integer $K$, that $$T_N(f)T_N(g)- \sum_{j=0}^K N^{-j}T_N(C_j(f,g))=O(N^{-K-1}).$$ The
functions $C_j$ are bilinear differential operators of degree less
than $2j$, and $C_0(f,g)=fg$. An explicit derivation of $C_j(f,g)$ is
given by the Proposition 6 of \cite{charles2003berezin}.

\subsection{A general localization result}
Using the $C^*$-algebra structure of Toeplitz operators, one can prove a fairly general
localization result:

\begin{prop}\label{prop:loc}
  Let $h$ be a smooth nonnegative function on $M$. Let $Z=\{h=0\}$, and suppose that $h$ vanishes exactly at order $2$ on $Z$, that is, there exists $c>0$ such that $h\geq c\dist(\cdot,Z)^2$. 
  
  Let $t>0$, and define $$V_t:=\{(m,v)\in X,\, \dist(m,Z)<t\}.$$
  
  For every $k\in \N$, there exists $C>0$ such that, for every $N\in \N$, for every $t>0$, and for every $u\in H_N(X)$ such that $T_N(h)u=\lambda u$ for some $\lambda \in \R$, one has
  
  $$\|u 1_{X\setminus V_t}\|^2_{L^2} \leq C\left(\cfrac{\max(\lambda,N^{-1})}{t^2}\right)^k\|u\|^2_{L^2}.$$
\end{prop}
\begin{rem}
  Here $M$ is a K\"ahler manifold, so $\dist$ is the Riemannian
  distance, but since $M$ is compact, the condition on $h$ does not
  depend on the chosen Riemannian structure.
\end{rem}

\begin{preuve}
	By a trivial induction, the $k$-th star power of a symbol $f$
        is of the form $$f^{\star k}=f^k+\eta C_{1,k}(f,\cdots, f) + \eta^2C_{2,k}(f,\cdots,f)+\ldots,$$ where $C_{i,k}$ is a $k$-multilinear differential operator of order at most $2i$.
	
	We want to study $C_{i,k}(h,\cdots, h)$ for $i\leq k$. The function $h$ is smooth and nonnegative, hence $\sqrt{h}$ is a Lipschitz function. In other terms, there exists $C$ such that, for every $(x,\xi)\in TM$ with $\|\xi\|\leq 1$, one has $|\partial_{\xi}h(x)|\leq C\sqrt{h(x)}.$
	In local coordinates, the function $C_{i,k}(h,\cdots, h)$ is a
        sum of terms of the form $a\partial^{\nu_1}h \partial^{\nu_2}h
        \ldots \partial^{\nu_k}h$, where $\sum_{j=1}^k|\nu_j|\leq 2i$ and $a$ is smooth.
	\begin{itemize}
		\item If $\nu_j=0$, then $\partial^{\nu_j}h =h$.
		\item If $|\nu_j|=1$, then $|\partial^{\nu_j}h| \leq C\sqrt{h}$.
		\item If $|\nu_j|\geq 2$, then $|\partial^{\nu_j}h| \leq C$.
	\end{itemize} 
	Hence $|a\partial^{\nu_1}h \partial^{\nu_2}h
        \ldots \partial^{\nu_k}h| \leq Ch^{k-\frac 12
          \sum_j\min(2,|\nu_j|)}$, moreover $\sum_j\min(2,|\nu_j|) \leq
      \sum_j |\nu_j|\leq 2i$, from which we can conclude:
	\begin{equation*}
		|C_{i,k}(h,\cdots, h)| \leq Ch^{k-i}.
	\end{equation*}
	
	This means that, for every $k\geq0$, the function $h^{\star k}$ is of the form:
	$$h^{\star k} = h^k+\sum_{i=1}^{k-1}\eta^{i}f_{i,k} + \eta^{-k}g(\eta),$$ where $g$ is bounded independently on $\eta$ and where, for each $i$ and $k$ there exists $C$ such that $|f_{i,k}| \leq Ch^{k-i}$.
	
	Using this, we can prove by induction on $k$ that there exists $C_k$ such that, for every $N$ and for every eigenvector $u$ of $T_N(h)$ with eigenvalue $\lambda$, one has $$|\langle u,h^ku\rangle| \leq C_k\max(\lambda,N^{-1})^k\|u\|^2.$$ Indeed, this is clearly true for $k=1$, because $\langle u,hu\rangle = \lambda\|u\|^2$.
	
	Let us suppose that, for all $1\leq i \leq k$, there exists $C$
        such that $$|\langle u,h^{k-i}u\rangle| \leq
        C\max(\lambda,N^{-1})^{k-i}\|u\|^2.$$ Because $u$ is an eigenvector for $T_N(h)$, it is an eigenvector for its powers, hence $$T_N(h^{\star k})u= T_N(h)^ku+O(N^{-\infty}) = \lambda^ku+O(N^{-\infty}).$$
	
	Replacing $h^{\star k}$ by its expansion and using the fact that $h\geq 0$, we find:
	$$|\langle u,h^ku\rangle| \leq \lambda^k\|u\|^2
        +\sum_{i=1}^{k-1}N^{-i}\langle u,f_{i,k}u\rangle +
        CN^{-k}\|u\|^2.$$ Here we used the fact that the function $g$
        is bounded.

	Now recall $|f_{i,k}|\leq C_{i,k}h^{k-i}$, and the induction hypothesis: $$|\langle u,h^{k-i}u\rangle| \leq C_i\max(\lambda,N^{-1})^{k-i}\|u\|^2$$ for every $i>0$. Hence
	$$|\langle u,h^ku\rangle| \leq C\max(\lambda,N^{-1})^k\|u\|^2 +\sum_{i=1}^{k-1}C_{i,k}C_iN^{-i}\max(\lambda,N^{-1})^{k-i}\|u\|^2,$$ hence there exists $C_k$ such that $|\langle u,h^{k}u\rangle |\leq C_k\max(\lambda,N^{-1})^k\|u\|^2$.
	
	Now we can conclude: for every $k$, there exists $C$ such that, for every $t>0$ one has $$\forall z\notin V_t,\,h^k\geq Ct^{2k}.$$
	
	Finally, for every $k$ there exists $C$ such that, for every $N\in \N$, $t>0$ and $u$ an eigenvector of $T_N(h)$ with eigenvalue $\lambda$, there holds 
	$$\|u 1_{X\setminus V_t}\|^2_{L^2} \leq C\left(\cfrac{\max(\lambda,N^{-1})}{t^2}\right)^k\|u\|^2_{L^2}.$$
	\end{preuve}

	Recalling Definition \ref{defn:precon}, let us specialize the Proposition \ref{prop:loc} to the case $\lambda=O(N^{-1})$ and $t=N^{-\delta}$ for $0<\delta<1/2$:
	\begin{corr}
		Let $u=(u_N)_{N\in \N}$ be a sequence of unit
                eigenvectors of $T_N(h)$, with sequence of eigenvalues
                $\lambda_N=O(N^{-1})$. If $h$ vanishes at order two on
                its zero set, then $u$ concentrates on this set.
	\end{corr}
	
	We can reformulate the Proposition \ref{prop:loc} in these
        terms: if $h$ is a positive smooth function on $M$, which
        vanishes at order two on its zero set, then any sequence of normalized eigenvectors of $T_N(h)$ whose eigenvalues are $O(N^{-1})$ concentrates on the zero set of $h$.
	
	\begin{rem}$\,$
	 \begin{itemize}
\item          An independent work by Charles and Polterovich, that
  appears partially in \cite{charles2015sharp}, treats the case of an
  eigenvalue close to a regular value of the symbol, with a result
  very similar to Proposition \ref{prop:loc}.
\item The proof of Proposition \ref{prop:loc} uses cancellation at
  order two only when dealing with $V_t$. Indeed, a more general result
  is
$$\|u1_{X\setminus V_t}\|^2_{L^2}\leq C
\left(\cfrac{\max(\lambda,N^{-1})}{\max(h(x),x\in
    V_t)}\right)^k\|u\|^2_{L^2},$$ which holds for any smooth $h$ and any eigenfunction $u$ of $T_N(h)$ with eigenvalue $\lambda$.
\end{itemize}
	\end{rem}
\subsection{Quadratic symbols on the Bargmann spaces}
Toeplitz operators can also be defined in the Bargmann spaces setting,
but one should be careful about the domain of such operators.

This section is devoted to a full survey of the quadratic case, which
is very useful as a model case for the general setting. Let $q$ be
a positive definite quadratic form on $\C^n$. Let $$\mathcal{A}_N=\left\{f\in \mathcal{B}_N,
\sqrt{q(\cdot)}f(\cdot)\in L^2(\C^n)\right\}.$$ Then
$\mathcal{A}_N$ is a dense subspace which contains the image of
$\mathcal{D}$ by the isomorphism between $\mathcal{B}_1$ and $\mathcal{B}_N$. It
is the domain of the positive quadratic form $t_N:(u,v)\mapsto \int
qu\overline{v}$, and $\mathcal{A}_N$ is closed for the norm
$\|u\|^2_{\mathcal{A}_N}=\|u\|^2_{L^2}+t_N(u,u)$. Moreover, the injection
$$(\mathcal{A}_N,\|\cdot\|_{\mathcal{A}_N})\hookrightarrow (\mathcal{B}_N,
\|\cdot\|_{L^2})$$ is compact. Using the usual results of spectral
theory, the associated operator $T_N^{flat}(q)$ is positive and has compact
resolvent. The spectrum of $T_N^{flat}(q)$ thus consists of a sequence of nonnegative
eigenvalues, whose only accumulation point is $+\infty$.

Observe that, since $q$ is even, $T_N^{flat}(q)$ sends even functions
to even functions, and odd functions to odd functions. Moreover, $q$ is
$2$-homogeneous. Recalling that the normalized scaling on $\C^n$ by a factor
$N^{1/2}$ sends $\mathcal{B}_N$ into $\mathcal{B}_1$, the conjugation
by this scaling
sends $T_N^{flat}(q)$ to $N^{-1}T^{flat}_1(q)$.

\begin{prop}
The first eigenvalue $\mu_N$ of $T_N^{flat}(q)$ is
  simple.
\end{prop}
\begin{preuve}
  As $q$ is positive a.e, the quadratic form $t_N$ is strictly convex, hence the first
  eigenvalue is simple.
\end{preuve}
We now compare Toeplitz quantization with Weyl quantization for
quadratic symbols.
Let $Op_{W}^{N^{-1}}(q)$ denote the Weyl quantization of $q$, as a symbol
  in $T^*\R^n\simeq \C^n$, with semiclassical parameter $N^{-1}$:
  \begin{equation*}
    Op_{W}^{N^{-1}}(q)u(x)=\cfrac{N^n}{(2\pi)^n}\int e^{iN\langle
      \xi,x-y\rangle}q
\left(\xi,\cfrac{x+y}{2}\right)u(y) \dd y \dd \xi.
  \end{equation*}
Recall that $B_N$ is the $N$-th Bargmann transform.

\begin{prop}\label{prop:quantequiv}
   $B_NT_N^{flat}(q)B_N^{-1}=Op^{N^{-1}}_{W}(q)+N^{-1}\cfrac{\tr(q)}{2}$.
  
  In particular, the first eigenvalue of $T_N^{flat}(q)$ is positive.
\end{prop}
\begin{preuve}
These computations belong to the folklore on the topic. Nevertheless,
for the comfort of the reader, we recover them explicitly.

It is sufficient to consider the $N=1$ case which is conjugated with the
general case through the usual scaling: indeed $Op^{N^{-1}}_W(q)=N^{-1}Op^1_W(q)$.

Here we shorten the notations for the momentum operators: on the Bargmann side, we let
$\mathfrak{d}_j=\partial_{z_j}+\frac 12 \overline{z_j}$; on the $\R^n$ side, we
let $\partial_{x_j}=\frac 1i \,\frac{\partial}{\partial x_j}$.

	Let $j,k$ be two indices in $[|1,n|]$. 
	
	If $q:z\mapsto z_jz_k=(x_j+iy_j)(x_k+iy_k)$, then $\tr(q)=0$, so the two operators should coincide. $T_1^{flat}(q)$ is the operator of multiplication by $z_jz_k$. This operator is conjugated via $B_1$ to the operator $(x_j+i\partial_{x_j})(x_k+i\partial_{x_k})=x_jx_k-\partial_{x_j}\partial_{x_k}+ix_j\partial_{x_k}+i\partial_{x_j}x_k$. Moreover, the Weyl quantization of $q$ is the operator $$Op^1_W(q) =x_jx_k-\partial_{x_j}\partial_{x_k} + \frac i2(\partial_{x_k}x_j+x_j\partial_{x_k}+\partial_{x_j}x_k+x_k\partial_{x_j}).$$ These two operators coincide whether $j=k$ or not.
	
	The case $q:z\mapsto \overline{z_j}\overline{z_k}=
        (x_j-iy_j)(x_k-iy_k)$ is the adjoint of the previous one.
	
	If $q:z\mapsto z_j\overline{z_k}=(x_j+iy_j)(x_k-iy_k)$, then $\tr(q)=2\delta^j_k$. In that case, $T_1^{flat}(q)=\mathfrak{d}_kz_j$. This operator is conjugated to $(x_k-i\partial_{x_k})(x_j+i\partial_{x_j})$. The Weyl quantization of $q$ is $$Op^1_W(q) =x_jx_k+\partial_{x_j}\partial_{x_k} + \frac i2(-\partial_{x_k}x_j-x_j\partial_{x_k}+\partial_{x_j}x_k+x_k\partial_{x_j}).$$ The two operators coincide when $k\neq j$, and when $k=j$ the difference is $1$.
	
	From the conjugation, it is clear that the first eigenvalue of
        $T_N^{flat}(q)$ is positive, because the Weyl quantization of $q$
        is nonnegative (see Remark \ref{rem:mu}) and $\tr(q)>0$. 
	\end{preuve}

Because $T_N^{flat}(q)$ is conjugated to $N^{-1}T_1^{flat}(q)$, one has  $\mu_N=N^{-1}\mu_1$, and for some $C>0$,
  \begin{equation*}
    \dist(\mu_N,Sp(T_N^{flat})\setminus \{\mu_N\})=CN^{-1}.
  \end{equation*}

The first eigenvalue $\mu_1$ of $T_1(q)$ depends on $q$, but is
invariant under an unitary change of variables on $\C^n$. From now on
we will use the notation $\mu(q)$ to denote $\mu_1$.

\begin{rem}\label{rem:mu}
  The computation of $\mu(q)$ is non-trivial. As explained in
  \cite{Melin1971}, Lemma 2.8, or as a direct consequence of the
  classification in \cite{williamson1936algebraic}, the first
  eigenvalue of $Op^1_W(q)$ can
  be obtained the following way: let $M\in S_{2n}^{+}(\R)$ denote the
  symmetric matrix associated with $q$ in the canonical
  coordinates. Let $J$ be the matrix of the symplectic structure: 
$$J=\begin{pmatrix}
    0&-Id\\Id&0
  \end{pmatrix}.$$

  Then the matrix $JM$ is skew-symmetric with respect to the scalar
  product given by $M$. Hence $A=iJM$ can be diagonalized; the eigenvalues of $A$ appear in pairs of opposite
  values $\lambda$ and $-\lambda$. Then $\mu$ is the sum of the positive diagonal
  elements of $A$.
In particular, this explicit formulation shows that $Op^{N^{-1}}_W(q)$
is nonnegative.
\end{rem}

We can use Proposition \ref{prop:quantequiv} to transpose well-known
results for the quantization of quadratic symbols to the Bargmann
case. Since $\mu(q)$ is simple, the operator $T^{flat}_1(q)-\mu(q)$
has a continuous inverse on the orthogonal set of the associated
eigenfunction. This inverse sends $\mathcal{D}$ into itself, because one
can build a Hilbert base of $\mathcal{D}$ with eigenfunctions of
$T^{flat}_1(q)$. Moreover the eigenfunction associated with $\mu(q)$ is
even.

\begin{rem}
  To illustrate the Proposition \ref{prop:quantequiv}, we solve
  completely the $n=1$ case. Up to a $U(1)$ change of variable, any
  real quadratic form on $\C$ can be written as $\alpha x^2+\beta
  y^2$. The associated Weyl operator is $-\beta \Delta + \alpha x^2$,
  with first eigenvalue $\sqrt{\alpha \beta}$. 
On the other hand, the first eigenfunction of
$\cfrac{\alpha-\beta}{4}(z^2+\mathfrak{d}^2) +
\cfrac{\alpha+\beta}{2}\,\mathfrak{d}z$ is a squeezed state of the form
$z\mapsto e^{-\frac 12 |z|^2}e^{\frac{\lambda}{2}z^2}$, with
$\lambda=\cfrac{(\sqrt{\alpha}-\sqrt{\beta})^2}{\alpha-\beta}$ (or
$\lambda=0$ in case $\alpha=\beta$). The associated eigenvalue is
then $\cfrac{(\sqrt{\alpha}+\sqrt{\beta})^2}{2}$. The difference is
$\cfrac{\alpha+\beta}{2}$, which is exactly half of the trace of $q$.
\end{rem}

\begin{rem}
  If instead of $T_N(h)$ one would consider $T_N(h-\frac{\Delta
  h}{2N})$, as in \cite{charles2003berezin}, then the Toeplitz quantization
  of a quadratic form would be exactly conjugated to its Weyl
  quantization: indeed $\tr(q)=\Delta q$. We recover in this
  particular case the computations in \cite{le2014singular}.
\end{rem}

\section{The first eigenvalue}

This section is devoted to the proof of Theorem A.

Let $P_0\in M$, one can find normal coordinates from a neighbourhood
of $P_0$ to a neighbourhood of $0$ in $\C^n$.
If at
$P_0$ a non-negative function $h$ vanishes with positive Hessian, the 2-jet of $h$ at $P_0$ maps to
a positive quadratic form $q$ on $\C^n$, up to a $U(n)$ change of
variables. Hence, the map associating to $P_0$ the first
eigenvalue $\mu$ of the model quadratic operator $T_N^{flat}(q)$ is
well-defined. From now on, we will also call $\mu$ this map.

The method of proof for Theorem A is then as follows: for each vanishing point
$P_0$, we construct a sequence of functions which concentrates on
$P_0$, consisting of almost eigenfunctions of $T_N(h)$, and for which the
associated sequence of eigenvalues is equivalent to
$N^{-1}\mu(P_0)$ as $N\to +\infty$. We then show a positivity estimate for
eigenfunctions concentrating on a single well. The uniqueness and the spectral gap
property follow from a similar argument. At every step, we compare
$T_N(h)$ with the operator on $\mathcal{B}_N$ whose symbol is the
Hessian of $h$ at the point of interest.

\subsection{Existence}
We let $h$ denote a smooth function satisfying the wells
condition. At every cancellation point of $h$, we will find a candidate for the ground state of $T_N(h)$. Instead of finding exact eigenfunctions, we search for approximate eigenfunctions. This is motivated by the following lemma:
\begin{lem}
	Let $T$ be a self-adjoint operator on a Hilbert space, let
        $\lambda \in \R$, and $u\in D(T)$ with norm $1$.
	
	Then $\dist(\lambda, Sp(T))\leq \|T(u)-\lambda u\|$.
\end{lem}

Let $P_0 \in M$ be a point where $h$ vanishes. Let $\rho$ be a local
map of normal coordinates in a neighbourhood of $\pi^{-1}(P_0)$. Let
$\Omega_N$ be the set of $z\in \C^n$ such that $(z/\sqrt N,0)$
belongs to the domain of $\rho$.
Recall from equation (\ref{eq:SNexpfixed}) that, for every $N\in \N$
and every $z,w \in \Omega_N$, there holds:
\begin{multline}\label{eq:SNexpscaled}
N^{-n}S_N^{P_0}\left(\frac{z}{\sqrt{N}},\frac{w}{\sqrt{N}}\right)\\=\Pi_1(z,w)\left(1+\sum_{k=1}^KN^{-k/2}b_k(z,w))\right)+R_K(z,w,N) + O(N^{-\infty}).
\end{multline}
Here the $b_j$'s are polynomials of the same parity as $j$, and there
exist $C>0,m>0$ such that, for every $(z,w,N)$ as above: $$|R_K(z,w,N)|\leq CN^{-(K+1)/2}e^{-C'|z-w|}(1+|z|^m+|w|^m).$$

The main proposition is
\begin{prop}
	\label{prop:well-ex}
	There exists a sequence $(u_j)_{j\geq 0}$ of elements of $\mathcal{S}(\C^n)$, with $\langle u_0,u_k\rangle = \delta^0_k$, and a sequence $(\lambda_j)_{j\geq 0}$ of real numbers, with $\lambda_0=\mu(P_0)$ and $\lambda_j=0$ for $j$ odd, such that, for each $K$ and $N$, if $u^K(N)\in L^2(X)$ and $\lambda^K(N)\in \R$ are defined as:
	\begin{align*}
		u^K(N)(\rho(z,\theta))&:=e^{iN\theta}N^{n}\sum_{j=0}^K N^{-j/2}u_j(\sqrt{N}z),\\
		u^K(N)&\text{ is supported in the image of $\rho$},\\
		\lambda^K(N)&=N^{-1}\sum_{j=0}^K N^{-j/2}\lambda_j,
	\end{align*}
	there holds, as $N
	\to +\infty$, $$\|S_NhS_Nu^K(N)-\lambda^K(N)u^K(N)\|_{L^2(X)}=O(N^{-(K+3)/2}).$$
\end{prop}

\begin{rem}
	The functions $u^K(N)$ do not lie inside $H_N(X)$, because they are identically zero on an open set. Nevertheless, the operator $S_NhS_N$ on $L^2(X)$ decomposes orthogonally into $T_N(h)$ on $H_N$, and $0$ on its orthogonal. Hence a nonzero eigenvalue of $S_NhS_N$ must correspond to an eigenvalue of $T_N(h)$ with same eigenspace. The same holds for almost eigenvalues.
	
	Introducing $\lambda^K$ as a polynomial in $N^{-1/2}$ whose odd terms vanish may seem surprising. However, in the proof, we construct $\lambda^K$ as a polynomial in $N^{-1/2}$, as we do for $u^K$. The fact that it is a polynomial in $N^{-1}$ is due to parity properties.
\end{rem}
\begin{preuve}
Let us solve the successive orders of $$(S_NhS_N-\lambda^K(N))u^K(N) \approx 0.$$

We write the Taylor expansion of $h$ around $P_0$ at order $K$ as $$h(x)=q(x)+\sum_{j = 3}^Kr_j(x)+E_K(x).$$ Because of equation (\ref{eq:SNexpscaled}), the kernel of $S_NhS_N$, read in the map $\rho$, is:
\begin{multline}
\label{eq:existlocalexp}
N^{-n}e^{i(\phi-\theta)}S_NhS_N\left(\rho\left(N^{-1/2}z,N^{-1}\theta\right),\rho\left(N^{-1/2}w,N^{-1}\phi\right)\right)\\ = N^{-1}\int  \left(q(y)+\sum_{k=1}^{K-2}N^{-k/2}r_{k+2}(y)+NE_K(N^{-1/2}y)\right)\\\times\left[\Pi_1(z,y)\left(1+\sum_{j=1}^{K} N^{-j/2}b_j(z,y)\right)+R_K(z,y,N)\right]\\\times\left[\Pi_1(y,w)\left(1+\sum_{l=1}^{K} N^{-l/2}b_{l}(y,w)\right)+R_K(y,w,N)\right]\dd y\\ + O(N^{-\infty}).
\end{multline}

Let us precisely write down the $K=0$ and $K=1$ case.

The dominant order (that is, $N^{-1}$) of the right-hand side is simply:
\begin{equation*}
(z,w)\mapsto N^{-1}\int_{\C^n} \Pi_1(z,y)q(y)\Pi_1(y,w) \dd y.
\end{equation*}

It is $N^{-1}$ times the kernel of the Toeplitz operator $Q=T^{flat}_1(q)$ on $B_1$ associated to the quadratic symbol $q$, which we studied in Subsection 3.3. Its resolvant is compact, the first eigenvalue $\mu(P_0)$ is simple, and if $u_0$ is an associated eigenvector, the operator $Q-\mu(P_0)$ has a continuous inverse on $u_0^{\perp}$ which sends $\mathcal{D}$ into itself. Moreover $u_0$ is an even function.

This determines $u_0$ and $\lambda_0=\mu(P_0)$. Here $u_0\in
\mathcal{D}$, so we can truncate the function $(z,\theta) \mapsto
e^{iN\theta}N^{n}u_0(N^{1/2}z)$ to a function supported on the domain
of $\rho$, with only $O(N^{-\infty})$ error. The push-forward by
$\rho$ of this truncation, extended by zero outside the image of
$\rho$, is denoted by $u^0(N)$.

Now $u_0 \in \mathcal{D}$ so $u^0$ concentrates on $P_0$. The error is thus:
\begin{multline*}
\|S_NhS_Nu^0(N)-N^{-1}\lambda_0u^0(N)\|_{L^2(X)}^2 \\\leq
CN^{-2}\hspace{-4pt}\int_{\Omega_N^3} A(z,y,w,N)^2|u_0(w)|^2 \dd y \dd w \dd z+O(N^{-\infty}),
\end{multline*}
where \begin{multline*}A(z,y,w,N)=N|E_2(N^{-1/2}y)\Pi_1(z,y)\Pi_1(y,w)|\\+h(y)\biggl(|R_0(z,y,N)|\Pi_1(y,w)+|R_0(y,w,N)|\Pi_1(z,y)\\+|R_0(z,y,N)R_0(y,w,N)|\biggl).\end{multline*}
Here, $E_2$ is a Taylor remainder of order $3$ on a compact set,
so $$|NE_2(N^{-1/2}y)|\leq C|y|^3N^{-1/2}.$$ Moreover, recall
that, on $\Omega_N^2$, one has $$|R_0(z,y,N)|\leq CN^{-1/2}e^{-C'|z-y|}(1+|z|^m+|y|^m).$$ Hence,
on $\Omega_N^3$, there holds:
$$|A(z,y,w,N)|\leq CN^{-1/2}e^{-C'|z-y|-C'|y-w|}(1+|z|^m+|y|^m+|w|^m).$$
Because $u_0\in \mathcal{D}$, one deduces:
\begin{align*}
&N^3\int_X |S_NhS_Nu^0-N^{-1}\lambda_0u^0|^2 \\&\leq C\int_{\Omega_N^3} e^{-2C'|z-y|-2C'|y-w|}(1+|z|^{2m}+|y|^{2m}+|w|^{2m})|u_0(w)|^2 \dd y\dd z\dd w\\&\hspace{25em} + O(N^{-\infty})\\
&\leq C\left(\int_{\C^n} |v|^{2m}e^{-C'|v|}\dd v\right)^2\int_{\C^n} |w|^{2m}|u_0(w)|^2 \dd w+O(N^{-\infty})\\&\leq C.
\end{align*}

This method (estimating an error kernel using polynomial
growth and off-diagonal exponential decay) will be used repeatedly again.

From there we deduce that $u_0$ is an approximate eigenvector:
\begin{equation*}
\|S_NhS_Nu_0(N)-N^{-1}\lambda_0u_0(N)\|_{L^2(X)}=O(N^{-3/2}).
\end{equation*}

This proves the proposition for the case $K=0$.

The construction of $u_1$ and $\lambda_1$ is different, moreover
there are supplementary error terms. The term of order $N^{-3/2}$ in the right-hand side of equation
(\ref{eq:existlocalexp}) is:
$$(z,w)\mapsto N^{-3/2}\int_{\C^n}
\Pi_1(z,y)[r_3(y)+q(y)(b_1(z,y)+b_1(y,w))]\Pi_1(y,w)\dd y .$$

Let $J_1$ denote the operator with kernel as above. We are trying to find $u_1$ and $\lambda_1$ such that
\begin{equation}\label{eq:existorder1}(Q-\lambda_0)u_1+ J_1u_0 = \lambda_1 u_0,\end{equation}
with the supplementary condition that $\langle u_1, u_0\rangle =0$: indeed if $(u_1,\lambda_1)$ is a solution of equation (\ref{eq:existorder1}), then so is
$(u_1+cu_0,\lambda_1)$ for any $c\in \C$. The orthogonality condition
makes the solution unique as we will see.

The functions $r_3$, $q$ and $b_1$ are polynomials, so
$J_1(\mathcal{D}) \subset \mathcal{S}(\C^n)$. This ensures that the problem
is well-posed. Note that $J_1$ does not map $\mathcal{D}$ into
holomorphic functions; this is because the normal map $\rho$ does not
preserve the holomorphic structure.

Now $r_3$ and $b_1$ are odd, so $J_1u_0$ is odd. In particular,
$\langle u_0, J_1u_0 \rangle=0$, and because $Q$ is self-adjoint,
$\langle u_0,(Q-\lambda_0)u_1 \rangle = 0$.
From this we deduce that $\lambda_1\|u_0\|^2=0$, hence $\lambda_1=0$.

To find $u_1$, we use again the fact that $J_1u_0$ is orthogonal to
$u_0$. Since $\lambda_0$ is a simple eigenvalue, $Q-\lambda_0$ is
invertible from $u_0^{\perp}$ to itself, and maps $\mathcal{S}\cap
u_0^{\perp}$ to itself. Hence there exists a unique
$u_1\in \mathcal{S}$ orthogonal to $u_0$, such that $(u_1,0)$ solves
(\ref{eq:existorder1}). Moreover $u_1$ is odd.

Now we investigate the error terms. With $u^1$ and $\lambda^1$ as in
the statement, let
$$f^1(N)=(S_NhS_N-\lambda^1(N))u^1(N).$$

As $u_0$ and $u_1$ belong to $\mathcal{S}$, the function $u^1$
concentrates on $P_0$, and so does $f^1$. Hence it is sufficient
to control $f^1$ near $P_0$. After a change of
variables, one has:
\begin{multline*}N^{-n}e^{-i\theta}f^1(N)(\rho(N^{-1/2}z,N^{-1}\theta))
  = N^{-2}J_1u_1(z)\\
 + \int \hspace{-0.2em}
  \Pi_1(z,y)\Pi_1(y,w)E_3(\frac{y}{\sqrt{N}})(1+\frac{b_1(z,y)}{\sqrt{N}})(1+\frac{b_1(y,w)}{\sqrt{N}})(u_0(w)+\cfrac{u_1(w)}{\sqrt{N}})
  \dd y \dd w\\
    + N^{-1}\hspace{-0.2em}\int\hspace{-0.2em}
    R_1(z,y,N)\Pi_1(y,w)(1+\frac{b_1(y,w)}{\sqrt{N}})(q(y)+\frac{r_3(y)}{\sqrt{N}})(u_0(w)+\cfrac{u_1(w)}{\sqrt{N}})
    \dd y \dd w\\ + N^{-1}\hspace{-0.2em}\int \hspace{-0.2em}
    \Pi_1(z,y)(1+\frac{b_1(z,y)}{\sqrt{N}})R_1(y,w,N)(q(y)+\frac{r_3(y)}{\sqrt{N}})(u_0(w)+\cfrac{u_1(w)}{\sqrt{N}})
    \dd y \dd w\\
+N^{-1}\hspace{-0.2em}\int \hspace{-0.2em}
R_1(z,y,N)R_1(y,w,N)(q(y)+\frac{r_3(y)}{\sqrt{N}})(u_0(w)+\frac{u_1(w)}{\sqrt{N}}) \dd y \dd w\\ + O(N^{-\infty}).
\end{multline*}

As $u_1\in \mathcal{S}$, the first line of the right-hand term is well-defined, and $\|N^{-2}J_1u_1\| = O(N^{-2})$.

There holds a uniform Taylor estimate on the domain of $\rho$:
$$E_3(y) \leq C |y|^4,$$ so $E_3(N^{-1/2}y)$ is bounded by $N^{-2}$
times a function with polynomial growth independent of $N$. In
particular, there exist $C,C',m>0$ such that, on $\Omega_N^3$:
\begin{multline*}|E_3(N^{-1/2}y)\Pi_1(z,y)\Pi_1(y,w)| \\ \leq
CN^{-2}e^{-C'|z-y|-C'|y-w|}(1+|z|^m+|y|^m+|w|^m).\end{multline*}

Of course the same type of estimate (with different $C$ and $m$)
applies if we multiply the left-hand side by $b_1(z,y)$, $b_1(y,w)$,
or both. Hence, following the last part of the $K=0$ case, we can
estimate the second line of the expansion of $f^1$ as $O_{L^2(X)}(N^{-2})$.

The three following lines are treated the same way: because $u_0$ and $u_1$
belong to $\mathcal{S}$, it is sufficient to prove estimates for the
error kernels, of the form $$|A(z,y,w,N)| \leq N^{-2}
Ce^{-C'|z-y|-C'|y-w|}(1+|z|^m+|y|^m+|w|^m),$$ which are easily checked.

We construct by induction on $K$ the following terms of the expansion.

For $j\geq 1$, we let $J_j:L^2(\C^n) \mapsto L^2(\C^n)$ the unbounded
and symmetric operator with kernel
\begin{equation*}
J_j(x,z)=\int_{\C^n}
\Pi_1(x,y)\Pi_1(y,z)\left(\sum_{\substack{k+l+m=j\\k,m,l\geq 0}}b_k(x,y)r_{2+l}(y)b_m(y,z)\right)\dd y.
\end{equation*}
Here we use the convention $b_0=1$, and $r_2=q$. The dense subspace
$\mathcal{S}(\C^n)$ is included in the domain of $J_j$, moreover
$J_j(\mathcal{S})\subset \mathcal{S}$ because all the $b_j$'s and
$r_l$'s are polynomials. Moreover $J_j$ has the same
parity as $j$.

Let $K\in \N$, and suppose we found functions $(u_k)_{k\leq K}\in \mathcal{S}$, orthogonal to $u_0$, and of the same parity as $k$, and real numbers $\lambda_k$ that vanish when $k$ is odd, and such that, for each $k\leq K$, there holds:
\begin{equation}
\label{eq:evexp}
(Q-\lambda_0)u_k+\sum_{j=1}^{k} J_ju_{k-j}= \lambda_ku_0 + \sum_{j=1}^{k-1}\lambda_ju_{k-j}.
\end{equation}

Let us find $u_{K+1}$, orthogonal to $u_0$, and $\lambda_{K+1}$ so
that equation (\ref{eq:evexp}) also holds for $k=K+1$.

Take the scalar product with $u_0$. As $Q$ is symmetric, the left-hand
side vanishes, and we get a linear equation in $\lambda_{K+1}$, whose
dominant coefficient is $\|u_0\|^2=1$. Hence $\lambda_{K+1}$ is
uniquely determined. Moreover, if $K+1$ is odd, then $J_ju_{K+1-j}$ and
$\lambda_ju_{K+1-j}$ are odd functions for every $j$, so their scalar products with
$u_0$ are zero, hence $\lambda_{K+1}=0$.

We now are able to find $u_{K+1}$ because we can invert $Q-\lambda_0$ on the orthogonal set of $u_0$. Finally, $u_{K+1}$ is of the same parity as $K+1$.

It remains to show that this sequence of functions $u$ corresponds to an approximate eigenvector of $S_NhS_N$.

Let $K\geq 0$, fixed in what follows. For each
$N\in \N$, we can build a function $u^K(N)$ on $X$, supported in the image of $\rho$ and such that, for $x$ in the image of $\rho$, one has
$u^K(N)(\rho(z,\theta))=e^{iN\theta}N^{n}\sum_{k=0}^KN^{-k/2}u_k(\sqrt{N}z)$. Note that $u^K(N)$ concentrates on $P_0$.

Let $$\lambda^K(N)=N^{-1}\sum_{k=0}^KN^{-k/2}\lambda_k.$$
We evaluate
$(S_NhS_N-\lambda^K(N))u^K(N)=:f^K(N)$. Consider an open set $V_1$, containing $P_0$, and compactly included in the image of $\rho$. One has $$\|f^K(N)\|_{L^{\infty}(^cV_1)}=O(N^{-\infty})$$ because $u^K(N)$ concentrates on
$P_0$.

To compute $f^K(N)$ in $V_1$,
we use the equation (\ref{eq:SNexpscaled}) at order $K$. A change of variables leads to:
\begin{multline*}
N^{-n}e^{-iN\theta}f^K(N)\left(\rho(N^{-1/2}z,\theta)\right)\\=N^{-1}\sum_{k=0}^KN^{-\frac k2}\left[(Q-\lambda_0)u_k(z)-\sum_{j=1}^k
J_ju_{k-j}(z)-\lambda_ku_0(z) -
\sum_{j=1}^{k-1}\lambda_ju_{k-j}(z)\right]\\
+N^{-1}\sum_{k=K+1}^{2K}N^{-\frac k2}\left[-\sum_{j=k-K}^{K}(J_j-\lambda_j)u_{k-j}(z)\right]\\
+\sum_{k,j,l=0}^KN^{-\frac{k+j+l}{2}}A_{j,l,N}u_k(z)
+\sum_{k,j=0}^KN^{-\frac{k+j}{2}}A'_{j,N}u_k(z)
+\sum_{k=0}^KN^{-\frac k2}A''_{N}u_k(z).
\end{multline*}

By construction, the first line of the right-hand term vanishes. The
second line is $O(N^{-(K+3)/2})$.
There are three error terms in the last line.
$A_{j,l,N}$ is the operator with kernel:
\begin{equation*}
A_{j,l,N}(z,w)=\int_{\Omega_N} \Pi_1(z,y)\Pi_1(y,w)b_j(z,y)b_l(y,w)E_K(N^{-1/2}y)\dd y .
\end{equation*}

The function $E_K$ is a Taylor remainder at order $K+3$, so there
exist constants $C>0,C'>0,m>0$ such that, on $\Omega_N^3$:
\begin{multline*}|\Pi_1(z,y)\Pi_1(y,w)b_j(z,y)b_l(y,w)E_K(N^{-1/2}y)|\\\leq CN^{-(K+3)/2}e^{-C'|z-y|+C'|y-w|}(1+|z|^m+|y|^m+|w|^m).\end{multline*}

Hence, for each function
$u\in \mathcal{S}$, one has $$\|A_{j,l,N}(u)\|_{L^2}=O(N^{-(K+3)/2}).$$ In particular it is true of the functions $u_k$.

$A'_{j,N}$ is the operator with kernel:
\begin{multline*}
A'_{j,N}(z,w)=\int_{\Omega_N} \Pi_1(z,y)b_j(z,y)h(N^{-1/2}y)R_K(y,w,N)\dd y\\+\int\Pi_1(y,w)b_j(y,w)R_K(z,y,N)h(N^{-1/2}y) \dd y .
\end{multline*}

One has $h(N^{-1/2}y) \leq CN^{-1}|y|^2$, so there are constants
$C>0,C'>0,m>0$ such that, on $\Omega_N^3$:
\begin{multline*}|\Pi_1(z,y)b_j(z,y)h(N^{-1/2}y)R_K(y,w,N)| \\\leq
CN^{-(K+3)/2}e^{-C'|z-y|-C'|y-w|}(1+|z|^m+|y|^m+|w|^m).\end{multline*}

As usual we get, for every $k$, that $$\|A'_{j,N}(u_k)\|_{L^2}=O(N^{-(K+3)/2}).$$

$A''_{N}$ is the operator with kernel
\begin{equation*}
A''_{N}(x,z)=\int_{\Omega_N^3}
R_K(x,y,N)h(N^{-1/2}y)R_K(y,z,N)\dd y.
\end{equation*}

Again there exist constants $C>0,C'>0,m>0$ such that, on $\Omega_N^3$:
\begin{multline*}|R_K(z,y,N)h(N^{-1/2}y)R_K(y,w,N)|\\\leq
CN^{-K-3}e^{-C'|z-y|-C'|y-w|}(1+|z|^m+|y|^m+|w|^m).\end{multline*}

To conclude, the $L^2$-norm of all the error terms is
$O(N^{-(K+3)/2})$.
\end{preuve}

From this proposition we conclude that, at every well $P$, \emph{there exists} an eigenvalue of $T_N(h)$ which has an asymptotic expansion in inverse powers of $N$, the dominant term being $N^{-1}\mu(P)$. In particular, the first eigenvalue of $T_N(h)$ is $O(N^{-1})$.
\subsection{Positivity}
The following proposition implies that the first eigenfunctions only concentrate on the wells that are minimal:
\begin{prop}
	\label{prop:well-pos}
  Let $(v_N)_{N\in \N}$ a sequence of normalized functions in $L^2(X)$. Suppose $v$ concentrates at a point $P_0$, on which $h$ vanishes. Then for each $\epsilon>0$ there
  exists $N_0$ and $C$ such that, if $N>N_0$, $$\langle v_N,hv_N \rangle \geq N^{-1}\mu(P_0)-CN^{-3/2+\epsilon}.$$
\end{prop}
\begin{preuve}
	Let $\delta< \frac 12$ be close to $\frac 12$. 
  Let $\rho$ denote a normal map around $P_0$. Then the
  sequence $(w_N)_{N>0}=(\rho^*v_N)_{N>0}$ is such that
  $\|w_N\|_{L^2(^cB(0,N^{-\delta}))}=O(N^{-\infty})$. Then one has as well:\begin{align*}\|\Pi_Nw_N\|_{L^2(^cB(0,2N^{-\delta}))}&=O(N^{-\infty})\\
  \|S_N^{P_0}w_N\|_{L^2(^cB(0,2N^{-\delta}))}&=O(N^{-\infty}).
  \end{align*}

 Using the Proposition \ref{prop:universal}, for $\delta$ close enough to $\frac 12$, if $\rho^*\Pi_N$ is a pull-back of $\Pi_N$ by $\rho$, one has $\|(S_N-\rho^*\Pi_N)v_N\|\leq CN^{-\frac 12+\epsilon}$. Hence, $$\|(S^{P_0}_N-\Pi_N)w_N\|\leq CN^{-\frac 12 + \epsilon}.$$

If $q$ is the Hessian
  of $h$ at $P_0$ read in the chosen coordinates, the spectrum of the model quadratic operator $\Pi_Nq\Pi_N$ is known: one has $$\langle
  w_N,\Pi_Nq\Pi_Nw_N\rangle \geq N^{-1}\mu(P_0)\|\Pi_Nw_N\|^2.$$ Moreover, on $B(0,2N^{-\delta})$ the following holds: $CN^{-2\delta}\geq h \geq q - CN^{-3\delta}$.
  
  Now, if $\delta$ is close enough to $\frac 12$, one has: \begin{align*}
  &\langle w_N,S_N^{P_0}hS_N^{P_0}w_N\rangle \\
  &\geq \langle w_N,S_N^{P_0}qS_N^{P_0}w_N\rangle-CN^{-3\delta}\\ 
  &= \langle w_N,S_N^{P_0}q\Pi_Nw_N\rangle + \langle w_N,S_N^{P_0}q(S_N^{P_0}-\Pi_N)w_N\rangle-CN^{-3\delta}\\
  &\geq \langle w_N,S_N^{P_0}q\Pi_Nw_N\rangle -CN^{-2\delta-\min(\delta,\frac 12 -\epsilon)}\\
  &= \langle w_N,\Pi_Nq\Pi_Nw_N\rangle +\langle w_N,(S_N^{P_0}-\Pi_N)q\Pi_Nw_N\rangle-CN^{-2\delta-\min(\delta,\frac 12 -\epsilon)}\\
  &\geq \langle w_N,\Pi_Nq\Pi_Nw_N\rangle -CN^{-2\delta-\min(\delta,\frac 12 -\epsilon)} \\ 
  &\geq N^{-1}\mu(P_0)-CN^{-2\delta-\min(\delta,\frac 12 -\epsilon)}.\end{align*}
Choosing $\delta$ such that $\delta \geq \frac 12 - \epsilon$
concludes the proof.
\end{preuve}

\begin{rem}
	In the proof, it appears that the condition of concentration
        on $P_0$ can be slightly relaxed. We only used the fact that,
        for some fixed $\delta$ determined by the geometry of $M$ and by $\epsilon$, one has $$\|v_N1_{\pi(x)\notin B(P_0,N^{-\delta})}\|_{L^2} = O(N^{-\infty}).$$ Thus, this proposition could be used in a more general context.
\end{rem}
\subsection{Uniqueness and spectral gap}
\begin{prop}
  \label{prop:well-un}
	Suppose $h$ satisfies the wells condition, and that there is
        only one well with minimal $\mu$. Then the approximate
        eigenvalues of proposition \ref{prop:well-ex} associated to
        this well correspond to the first eigenvalue $\lambda_N$ of
        $T_N(h)$, namely, for every $K\in \N$, there holds: $$|\lambda^K(N)-\lambda_N|=O(N^{-(K+3)/2}).$$ This eigenvalue is simple; moreover there exists $C>0$ such that, for $N$ large enough: $$\dist(\lambda_N, \Sp(T_N(h))\setminus \{\lambda_N\})\geq CN^{-1}.$$
\end{prop}
\begin{preuve}
	The proposition is equivalent to the claim that there exists
        $K$ such that the following is true: let $u_K(N)$ denote the approximate eigenvector of order $K$ associated to the well with minimal $\mu$. Let $F_N$ be the orthogonal complement of $u_K(N)$ in $H_N(X)$, and $p_N$ be the orthogonal projection from $H_N(X)$ to $F_N$. Then the operator $T_N^{\sharp}(h):F_N\to F_N$, defined as $T_N^{\sharp}(h)=p_NT_N(h)$, is bounded from below by $\lambda_N+CN^{-1}$.
	
	Let $v_N$ be a sequence of normalized eigenvectors of $T_N^{\sharp}(h)$, and $\mu_N$ the sequence of associated eigenvalues. One has $T_N(h)v_N=\mu_Nv_N+C_Nu_K(N)$. Because $u_K(N)$ is a sequence of normalized functions and $S_N$ is bounded, the sequence $C_N$ is bounded.
	
	Assume $\mu_N=O(N^{-1})$. In this slightly different setting, we
        can adapt the proof of Proposition \ref{prop:loc} using
        the fact that $u_K(N)$ is itself an almost eigenfunction of $T_N(h)$. There holds:
	$$T_N(h^{\star
          k})v_N=\mu_N^kv_N+C_N\sum_{j=1}^k\mu_N^{j-1}\lambda_N^{k-j}u_K(N)
        + O(N^{-(K+3)/2}).$$

We can proceed as in \ref{prop:loc} but the induction process stops
at $k=\frac{K+3}{2}$. One concludes that, for every $\epsilon>0$,
the $L^2$ norm of $v_N$ is $O(N^{-\frac{K+3
 -\epsilon}{4}})$ outside the union of balls centred at the
vanishing points of $h$,
and of radius $N^{-\frac 12 + \frac{\epsilon}{K+3}}$.
	
	In particular, if $P_0,P_1,\ldots, P_d$ denote the vanishing points of $h$, and
        $P_0$ is the only one with minimal $\mu$, one can decompose
        $v_N=v_{0,N}+v_{1,N}+\ldots + v_{d,N} + O(N^{-(K+3-\epsilon)/4})$,
        where each sequence $v_{i,N}$ concentrates on $P_i$. The
        proposition \ref{prop:well-pos} gives estimates for $v_{i,N}$
        if $i \neq 0$. Indeed $\mu(P_i)>\mu(P_0)$ by construction, and
        $\lambda_N \leq N^{-1}\mu(P_0) + O(N^{-3/2})$, so one can find
        $C>0$ small enough such that $N\lambda_N+C<\mu(P_i)$ for all $i\neq 0$ and
        for $N$
        large enough. Then $$\langle v_{i,N},S_NhS_Nv_{i,N}\rangle \geq (\lambda_N+ CN^{-1})\|v_{i,N}\|_2^2.$$
	
	Recall that $u_K(N)$ has an asymptotic expansion whose first
        term $u_0$ is the pull-backed ground state of the operator on
        the Bargmann space with quadratic symbol. This operator has a
        (fixed) nonzero specral gap. Moreover $\langle
        v_{0,N},u_K(N)\rangle = O(N^{-(K+3-\epsilon)/4})$ because $v_N$ is
        orthogonal to $u_K(N)$ and $u_K(N)$ concentrates only at
        $P_0$. Then for $C$ strictly smaller than the spectral gap of
        the quadratic operator $T_1^{flat}(q_0)$ at $P_0$, one has for $N$ large $$\langle v_{0,N},S_NhS_Nv_{0,N}\rangle \geq (\lambda_N+CN^{-1})\|v_{0,N}\|_2^2.$$
	
	The functions $v_{i,N}$ have disjoint supports, so that
        $\langle v_{i,N},S_NhS_Nv_{j,N}\rangle = O(N^{-\infty})$
        whenever $i \neq j$, and $\|v_N\|_2^2=\sum_j
        \|v_{j,N}\|_2^2+O(N^{-(K+3-\epsilon)/4})$. Thus the two
        inequalities allow us to conclude when $K\geq 2$.
\end{preuve}

\subsection{End of the proof}
It remains to show that, in the case where only one well $P_0$ has minimal $\mu$, then the ground state is $O(N^{-\infty})$ in a fixed neighbourhood of the other wells.
Let $K\in \N$. We have constructed in Subsection 4.1 a sequence
$(u_K(N))_{N\in \N}$ which vanishes outside a fixed neighbourhood of
$P_0$, and which is a sequence of approximate unit eigenvectors of
$T_N(h)$, with approximate eigenvalue $\lambda_K(N)$. One
has $$\lambda_K(N)= N^{-1}\mu(P_0)+ O(N^{-3/2}),$$
and $$\dist(\lambda_K(N),\Sp(T_N(h)) =O(N^{-(K+3)/2}).$$ Moreover we
proved in Subsection 4.3 that there can be only one eigenvalue of
$T_N(h)$ in $[0,N^{-1}(\mu(P_0)+C)]$ for some $C$, and that this
eigenvalue is simple. Hence, denoting $\lambda_{\infty}(N)$ this
sequence of eigenvalues, one has $$\lambda_{\infty}(N)=\min
\Sp(T_N(h)),$$
and $$|\lambda_{\infty}(N)-\lambda_K(N)|=O(N^{-(K+3)/2}).$$ Let
$U_{\infty}(N)$ denote a sequence of unit eigenvectors associated to
$\lambda_{\infty}(N)$, and decompose $u_K(N)= c(N)U_{\infty}(N)+
w_K(N)$, where $w_K(N)\perp
U_{\infty}(N)$. Then $$(T_N(h)-\lambda_{\infty}(N))w_K(N)=O(N^{-(K+3)/2}).$$
The operator $T_N(h)-\lambda_{\infty}(N)$ is invertible on
$U_{\infty}(N)^{\perp}$ and its inverse has a norm bounded by $N$, so
$w_K(N)=O(N^{-(K+1)/2})$. Since both $u_K(N)$ and $U_{\infty}(N)$ are normalized, one has $c(N) \to 1$.

Finally, if $V$ is a neighbourhood of another well, then $u_K(N)$ is zero on $V$, so that $$\|U_{\infty}(N)\|_{L^2(V)}=\|w_K(N)\|_{L^2(V)}=O(N^{-(K+1)/2}).$$ This concludes the proof.

\section{Eigenvalues in a scaled window}

This section is devoted to the proof of Theorem B. The method of
proof is very similar to that of Theorem A: we will exhibit
approximate eigenvectors, then show that they cover the low-energy spectrum.
\subsection{Approximate eigenvectors}
In the proof of the Proposition \ref{prop:well-ex}, the first guess for an approximate eigenvector of $T_N(h)$ was the first eigenvector of the model quadratic operator at one of the wells.
If, instead of the first eigenvector, we start from any eigenvector of the model quadratic operator, we can proceed the same way; however the recursion stops after one step, in general.

\begin{prop}\label{prop:ex}
	Let $P\in M$ on which $h$ cancels, and $Q$ be a model
        quadratic operator in some normal map $\rho$ around $P$. Let
        $\lambda$ be an eigenvalue of $Q$ and $E_{\lambda}$ the
        corresponding eigenspace.
	Then one can find a suitable orthonormal basis $(v_1,\ldots,
        v_d)$ of $E_{\lambda}$, functions
        $(w_1,\ldots, w_d)$ in $\mathcal{S}(\C^n)$ and real numbers $(b_1,\ldots,b_d)$ such that, for any integer $i
        \in [1,d]$, the function $$\tilde{v}_i(N):\rho(z,\theta)\mapsto
        N^ne^{iN\theta}(v_i(N^{1/2}z)+N^{-1/2}w_i(N^{1/2}z))$$ is such
        that 
$$S_NhS_N \tilde{v}_i(N)= N^{-1}\lambda + N^{-3/2}b_i + O(N^{-2}),$$
	
	Moreover, if $\dim E_{\lambda}=1$, then if $u_0$ is an
        eigenvector of $Q$, one can find a sequence of Schwartz
        functions $(u_k)_{k\geq 1}$, orthogonal to $u_0$, and a
        sequence of real numbers $(\lambda_k)_{k\geq 1}$, such that,
        for every $K>0$, the function $$u_K(N):\rho(z,\theta)\mapsto
        N^ne^{iN\theta}\sum_{k=0}^KN^{-k/2}u_k(N^{1/2}z)$$ is such
        that $$S_NhS_N u_K(N)=N^{-1}\lambda +
        N^{-1}\sum_{k=1}^{K/2}N^{-k}\lambda_k = O(N^{-(K+3)/2}).$$
\end{prop}

\begin{preuve}
Recall from Proposition \ref{prop:well-ex} that one can find an approximate eigenvector at any order, starting from the ground state $u_0$ of $Q$.

Let now $u_0$ denote an arbitrary eigenfunction of $Q$, which still
belongs to $\mathcal{D}$. Let $\lambda$ be the associated
eigenvalue. When $\lambda$ is simple, the operator
$Q-\lambda$ has a continuous inverse on $u_0^{\perp}$, so one can solve equation (\ref{eq:evexp}) at any order. Observe that $u_0$ is either even or odd, so that only negative integer powers of $N$ remain in the expansion of the eigenvalue.

If $Q-\lambda$ is not invertible on $u_0^{\perp}$, the equation
(\ref{eq:evexp}) can still be solved for $K=1$ if $u_0$ is one of the
vectors of a convenient basis of $E_{\lambda}$; but the construction
fails at higher orders. Consider an orthonormal basis $(v_1,\ldots, v_L)$ of the eigenspace $E_{\lambda}$. Suppose $u_0=v_l$. The equation (\ref{eq:evexp}) reads:
\begin{equation*}
	(Q-\lambda)u_1+J_1u_0=\lambda_1u_0.
\end{equation*}
Taking the scalar product with $u_0$ yields $\lambda_1(l)=\langle
v_l,J_1v_l\rangle$. But we also need to check that $0=\langle
v_l,J_1v_j\rangle$ for $l\neq j$. This is done by choosing an
orthogonal basis in which the corestriction of $J_1$ on $E_{\lambda}$
is diagonal (recall $J_1$ is symmetric and $E_{\lambda}$ is finite-dimensional). One can then find $u_1(l)$ in
$E_{\lambda}^{\perp}$. The proof of the error estimate is the same. To
conclude we let $b_l=\lambda_1(l)$ and $w_l=u_1(l)$.

Once the $K=1$ step is done, the basis $(v_1,\ldots, v_L)$ is
fixed. Let us try to solve equation (\ref{eq:evexp}) with $u_0=v_1$, for $K=2$. We
write
$$(Q-\lambda)u_2+J_2u_0+J_1u_1=\lambda_2u_0+\lambda_1u_1.$$
Taking the scalar product with $u_0$ yields $\lambda_2$ as previously:
$$\lambda_2=\langle u_0,J_2u_0\rangle + \langle u_0,J_1u_1\rangle.$$
Now recall $u_1$ is orthogonal to $E_{\lambda}$. If $v$ denotes an
element of $E_{\lambda}$ orthogonal to $u_0$, then one must
check $$\langle v,J_2u_0\rangle + \langle v,J_1u_1\rangle=0.$$ This
equation does not hold in general, hence the obstruction.
\end{preuve}

\subsection{Uniqueness}
Let $C'>0$, and $N\in \N$. Consider the set $e_N$ of approximate eigenvectors in
Proposition \ref{prop:ex}, such that $\lambda<C'$. Then
$E_N=span(e_N)$ is a subspace of $L^2(X)$, with
small energy: there exists $C_1$ such that, for every $N$, $$\max\{\langle u,T_N(h) u \rangle, u\in E_N,
  \|u\|^2_2=1\}<C'N^{-1}+C_1N^{-\frac 32}.$$

We claim that, reciprocally, any function approximately orthogonal
with $E_N$ has an energy bounded from below:
\begin{prop}
  Let $C'>0$. There exists $\epsilon_0>0$ and a function $\epsilon
  \mapsto N_0(\epsilon)$ such that, for
  $0<\epsilon<\epsilon_0$, the following is true. Let $v_N$ be a normalized eigenfunction of $T_N(h)$, with associated
  eigenvalue $\lambda_N$, and suppose that the angle
  between $v_N$ and $E_N$ is greater than $\cos^{-1}(\epsilon)$, that is, for
  every $u\in E_N$ normalized, one has $|\langle u,v_N
  \rangle|<\epsilon$.
Then for $N\geq N_0(\epsilon)$, one has $$\lambda_N \geq (C'-\epsilon)N^{-1}.$$

\end{prop}
\begin{preuve}

Let $P_0,\ldots, P_d$ denote the points at which $h$ cancels. If
$\lambda_N = O(N^{-1})$, then $v_N$ concentrates on the $P_i$'s. We
decompose $v_N=v_{0,N}+v_{1,N}+\ldots+v_{d,N}+O(N^{-\infty})$, where
each $v_{i,N}$ concentrates only on $P_i$.

Let $\rho_i$ be a normal map associated with $P_i$, and $q_i$ the
Hessian of $h$ at $P_i$ read in the map $\rho_i$. Let $E_{i,N}$ be
the span of eigenfunctions of $T_N^{flat}(q_i)$ whose eigenvalues are
less than $C'N^{-1}$. Then for $N$ large, for every normalized $u\in
E_{i,N}$, one has $|\langle\rho_i^* v_{i,N},u\rangle|\leq
2\epsilon$. Indeed functions in $E_{N}$ are $N^{-1/2}$-close to sums of
pull-backs of functions in $E_{i,N}$.

Hence, for $N$ large enough, $$\langle
\rho_i^*v_{i,N},\Pi_N(q_i-C'N^{-1})\Pi_N\rho_i^*v_{i,N}\rangle \geq
-C'N^{-1}(4\epsilon^2).$$

Since $v_{i,N}$ concentrates on $P_i$, one can deduce that, for $N$
large enough,
$$\langle v_{i,N},S_NhS_N v_{i,N}\rangle \geq
C'N^{-1}\|v_{i,N}\|^2-C'N^{-1}(5\epsilon^2),$$
hence
$$\langle v_N,S_NhS_Nv_N\rangle \geq
C'N^{-1}-C'N^{-1}(5(d+2)\epsilon^2).$$

To conclude, we let $\epsilon_0=\cfrac{1}{5(d+2)C'}$. Then for every
$\epsilon<\epsilon_0$, for $N$ large enough,
$$\langle v_N,S_NhS_Nv_N\rangle \geq
(C'-\epsilon)N^{-1}.$$
\end{preuve}

To conclude the proof of Theorem $B$, if the rank of the spectral
projector of $T_N(h)$ with interval $[0,CN^{-1}]$ was greater than
$K$, then one could find an eigenfunction of $T_N(h)$ which forms an
angle greater than $\cos^{-1}(N^{-1})$ with $E_N$, and with eigenvalue
less than $CN^{-1}$. This is absurd since $C<C'$.

\section{Acknowledgements}
The author thanks N. Anantharaman and L. Charles for their help and
encouragement while writing this article.

\appendix

\section*{Appendix : a proof for the off-diagonal estimate}
\setcounter{section}{1}
This last section is an appendix about Proposition
\ref{prop:ma4.18}. As we already explained, the knowledge of the
result is sufficient for our needs. However, as this proposition
appears in \cite{ma2006}, it is stated in a case that is much more
general than prequantum bundles on K\"ahler manifolds.

In this specific setting, and with a more direct approach, we propose to show a different version of this estimate, with a
somewhat stronger estimate on the remainder (see Proposition
\ref{prop:rem}). We also replace the normal maps of Definition \ref{defn:normap}
with Heisenberg maps, satisfying different assumptions. This version could be of use in
situations where it is crucial that the local map is a biholomorphism.

The proof relies on the theory of Fourier Integral Operators with
complex-valued phase functions, in the sense of H\"ormander
(\cite{hormander2003analysis}, section 7.8). Indeed, we will follow the lines of \cite{Shiffman2002} (restricting
ourselves to exact K\"ahler structures), which gives asymptotics at a
shrinking scale; we modify the proof in order to estimate the
remainder at a fixed scale, recovering results from \cite{charles2003berezin,Berman2008}.

The starting point in \cite{Shiffman2002} is the study by Boutet de
Monvel and Sj\"ostrand \cite{BoutetdeMonvel1975} of the
general Szeg\H{o} projector (Definition \ref{defn:Hardy-Szego}). The structure of the Szeg\H{o}
projector, for the boundary of a relatively compact open set, has been subject to
a thorough study (\cite{kohn1963harmonic, kohn1964harmonic,
  kohn1965extension, BoutetdeMonvel1974-1975, BoutetdeMonvel1975, BoutetdeMonvel1981}). Under the assumption of \emph{strong pseudo-convexity},
which is verified for the unit ball $D$ of $L^*$, the boundary of $D$
inherits a Riemannian metric from the Levi form (which is identical to
the one we use in this paper).
The projector $S$ is then a Fourier Integral Operator with complex phase, in the sense of H\"ormander \cite{hormander2003analysis}:
\begin{prop}[\cite{BoutetdeMonvel1975}]\label{prop:Boutet}Let $Y$ be
  the boundary of a strongly pseudo-convex, relatively compact open set in a complex
  manifold. Then there exists a skew-symmetric almost holomorphic complex phase function $\psi \in
  C^{\infty}(Y\times Y)$ (in the sense of \cite{hormander2003analysis}), whose critical set is $\diag(Y)$, and a
  classical symbol $$s\sim \sum_i t^{n-i}s_i\in C^{\infty}(Y\times Y
  \times \R^+),$$ such that the Schwartz kernel of the Szeg\H{o} projector on $Y$
  is $$S(x,y)= \int_0^{+\infty} e^{it\psi(x,y)}s(x,y,t)\dd t+ E(x,y),$$ where the function $E$ is smooth. Moreover the
  principal symbol $s_0$ is such that $s_0^2=h_{\psi}^{-1}$, where
  $h_{\psi}(x,y)$ is the Hessian of the function $$Y\times \R^+ \ni (z,\sigma)\mapsto \psi(x,z)+\sigma \psi(z,y)$$ at the critical point (which is unique and lies in a complex extension of $Y\times \R^+$).
\end{prop}
In this setting, ``almost holomorphic'' means that, near the diagonal $z=w\in Y$, one has $\overline{\partial}_z\psi(z,w)=O(|z-w|^{\infty})$.
The fact that the function $(z,\sigma)\mapsto \psi(x,z)+\sigma
\psi(z,y)$ has exactly one critical point in the complex extension of
$Y\times \R^+$, with nondegenerate Hessian, is encoded in the
requirements on $\psi$ to be a complex phase function in the sense of
H\"ormander.

In the specific case where $X$ is a circle bundle over $M$, one can
use the microlocal information on $S$ to deduce the asymptotics of its
Fourier components $S_N$. Indeed, the $N$-th Fourier component of a smooth function on a compact set has a sup norm bounded by $O(N^{-\infty})$. Thus, one has
$$S_N(x,y)= \iint
\exp(it\psi(x,r_{\eta}y)+iN\eta)s(x,r_{\eta}y,t)\dd t \dd \eta +
E_N(x,y),$$ where $\|E_N\|_{L^{\infty}}=O(N^{-\infty}).$
Here, as in the introduction, $r_{\eta}$ denotes the circle action on $X$.

As announced, we will deal with a less restrictive class of local
maps than the normal maps of Definition \ref{defn:normap}. Because we are dealing with exact K\"ahler manifolds, as opposed to the more general almost complex structure, we slightly
modify the definition of \cite{Shiffman2002}:

\begin{defn}
  Let $P_0\in M$. Let $U$ be a neighbourhood of $0$ in $\C^n$ and $V$
  be a neighbourhood of $P_0$ in $M$.

  A smooth diffeomorphism $\rho:U\times \R\to \pi^{-1}(V)$ is said to be an
  \emph{Heisenberg map} or map of \emph{Heisenberg coordinates} under
  the following conditions:
  \begin{itemize}
  \item $\pi(\rho(0,0))=P_0$;
  \item $\rho^*\omega(P_0)= \omega_0(0)$.
  \item $\overline{\partial}\rho=0$.
  \item $\rho(m,\theta)=r_{\theta}\rho(m,0)$.
  \end{itemize}
\end{defn}

 The crucial point is that, in these coordinates, the phase $\psi
$ from the Boutet-Sj\"ostrand theorem
reads, for all $(z,\theta)$ and $(w,\phi)$ in the domain of $\rho$
(cf. \cite{Shiffman2002}, equation 61): $$\psi(\rho(z,\theta),\rho(w,\phi))=i\left[1-A(z,w)e^{i(\theta-\phi)}\right],$$
  
  Here, the $2$-jet of $A$ is known at the origin
  (\cite{Shiffman2002}, Lemma 2.4): $$A(z,w)=1 - \frac 12 |z-w|^2 + i\Im(z\cdot \overline{w})+O(|z|^3,|w|^3).$$
  
  We will need to control the off-diagonal behaviour of
  $A$. Recall $$\Pi_1:(z,w)\mapsto \cfrac{1}{\pi^n}\exp\left(-\frac 12
    |z-w|^2+i\Im(z\cdot \overline{w})\right).$$ Up to a reduction of
  the definition set of $\rho$, the usual logarithm is well-defined,
  and we can define $R_A$ as the unique function such that
  $A/\Pi_1=\pi^n e^{R_A}$.
  \begin{prop}
  	\label{prop:remainder}
  	The two following estimates hold as $z,w\to 0$:
  	\begin{align*}
  		\Re(R_A)(z,w) &= O\left(|z-w|^2(|z|+|w|)\right)\\
  		\Im(R_A)(z,w) &= O\left(|z-w|(|z|^2+|w|^2)\right).
  	\end{align*}
  	In particular, up to a restriction of the Heisenberg map $\rho$ to a smaller neighbourhood of $P_0$, one has, for every $z$ and $w$ in the domain of $\rho$:
  	\begin{equation}
  	|A/\Pi_1|(z,w) \leq \pi^n e^{\frac 14|z-w|^2}.
  	\end{equation}
  \end{prop}
  \begin{preuve}
  	The functions $A$ and $\pi^n\Pi_1$ are equal up to order $2$ at $P_0$, so that $R_A(z,w)=O(|z|^3,|w|^3)$.
  	
  The two functions $A$ and $\pi^n\Pi_1$ are both smooth and are equal to
  $1$ on the diagonal. Moreover the first derivatives of both $\Re(A)$
  and $\Re(\Pi_1)$ vanish on the diagonal. For $\Pi_1$ this is a
  straightforward computation. For $A$ it comes from the fact that
  $\psi$ is a complex phase function whose critical set is the
  diagonal. It is also a natural consequence of the fact that
  $\partial_1A(z,z)=-\frac 12 \partial \phi(z)$ and
  $\overline{\partial}_1A(z,z)=\frac 12 \partial \phi(z)$, where
  $\phi$ is a complex potential: $i\partial \overline{\partial}\phi=\omega$. Hence there is a constant $C$ such that, for every $z$ and $w$ in the domain of $\rho$, there holds:
  \begin{align*}
  	|\Im(A-\pi^n\Pi_1)(z,w)| &\leq C|z-w|(|z|^2+|w|^2)\\
  	|\Re(A-\pi^n\Pi_1)(z,w)| &\leq C|z-w|^2(|z|+|w|).
  \end{align*}
  From which we deduce that
  \begin{align*}
	|\Re((A-\pi^n\Pi_1)^2)(z,w)| & \leq C|z-w|^2(|z|+|w|)\\
  	|\Im((A-\pi^n\Pi_1)^2)(z,w)| & \leq C|z-w|^3\\
  	|A-\pi^n\Pi_1|^3 \leq |z-w|^3.
  \end{align*}
  Now
  $$R_A=\log(A/\pi^n\Pi_1)=\cfrac{A-\pi^n\Pi_1}{\pi^n\Pi_1}-\cfrac 12 \, \left( \cfrac{A-\pi^n\Pi_1}{\pi^n\Pi_1}\right)^2+O\left(\left( \cfrac{A-\pi^n\Pi_1}{\pi^n\Pi_1}\right)^3\right).$$
  Taking the real and imaginary part of this equation, one deduces
  \begin{align*}
  \Re(R_A)(z,w) &= O\left(|z-w|^2(|z|+|w|)\right)\\
  \Im(R_A)(z,w) &= O\left(|z-w|(|z|^2+|w|^2)\right).
  \end{align*}
  In particular,
  \begin{equation*}
	  |A/\Pi_1|(z,w) \leq \pi^n e^{C|z-w|^2(|z|+|w|)}.
  \end{equation*}
  Reducing the domain of the Heisenberg map $\rho$ to a smaller neighbourhood of $P_0$, one gets, for every $z$ and $w$ in the domain of $\rho$:
  \begin{equation*}
  |A/\Pi_1|(z,w) \leq \pi^n e^{\frac 14|z-w|^2}.
  \end{equation*}
  \end{preuve}
In fact, the symbol $s$ of the operator can also be chosen to be very
simple in the given coordinates:
\begin{prop}
  In Heisenberg coordinates, the symbol $s$ of $S$ in proposition \ref{prop:Boutet} can be chosen to be factorized as: $$s(\rho(z,\theta),\rho(w,\phi),t)=e^{-in
    (\theta-\phi)}\xi(z,w,t),$$ where $$\xi(z,w,t) \sim \sum_{k=0}^{+\infty}t^{n-k}\xi_k(z,w)$$ and where each $\xi_k$ is a smooth function. Moreover the principal symbol $\xi_0$ does not vanish in
  a neighbourhood of $\diag(M)$.
\end{prop}
\begin{preuve}
  The expression of the phase $\psi$ in local coordinates gives
  immediatly that any derivative of order $\geq 2$ of the function
  $(t,z,\theta,w,\phi)\mapsto t\psi(\rho(z,\theta),\rho(w,\phi))$ is
  of the form $e^{i(\theta-\phi)}f(z,w,t)$ where $f$ is constant or
  linear wrt $t$. It follows that $h_{\psi}(\rho(z,\theta),\rho(w,\phi))= e^{2in
    (\theta-\phi)}g(z,w)$ for some function $g$. Hence, we can write $s_0(\rho(z,\theta),
\rho(w,\phi))=e^{-in
    (\theta-\phi)}\xi_0(z,w)$ for some smooth function $\xi_0$. Of course, any partial derivative
  of $s_0$ is also, in local coordinates, of the form $e^{-in(\theta-\phi)}f(z,w)$ for some function $f$.

  Let us assume that for $k\leq K$, each function $s_k$ reads in local
  coordinates as $e^{in(\theta-\phi)}\xi_k(z,w)$ for some smooth function $\xi_k$. The coefficient $s_{K+1}$ can be derived from
  $(s_i)_{i\leq K}$ via a
  stationary phase lemma, in which the differential operators come
  from the Taylor expansion of $\psi$. Thus, $s_{K+1}$ is a priori of the
  form $$s_{K+1}(\rho(z,\theta),\rho(w,\phi))=e^{-in(\theta-\phi)}\left(\sum_{j=-C}^Ce^{ik(\theta-\phi)}\xi_{K+1,j}(z,w)\right),
  $$
  where $C$ is finite (but depends on $K$) and the $\xi_{K+1,j}$ are smooth functions. 
  
  We can get rid of all coefficients except $j=0$
  by adding a convenient multiple of $\psi$. Indeed, the operator with symbol $(f+\psi g)t^k$ is equal, after integration by parts, to the operator with symbol $ft^k+ikgt^{k-1}$, modulo a smoothing operator.  
  For instance, replacing
  $s_{K+1}$ with $s_{K+1}+e^{-i(\theta-\phi)}\xi_{K+1,1}a(z,w)\psi$
  eliminates the $j=1$ term. 

We conclude by induction.
\end{preuve}

The $N$-th Fourier component $S_N$ of the Szeg\H{o} projector at a
point $(x,y)$ reads $$S_N(x,y)=\iint
\exp(it\psi(x,r_{\eta}y)+iN\eta)s(x,r_{\eta}y,t)\dd t \dd \eta+ O(N^{-\infty}).$$ A change of
variables leads to $$S_N(x,y)=N\iint
\exp(iN(t\psi(x,r_{\eta}y)+\eta))s(x,r_{\eta}y,Nt) \dd t \dd \eta + O(N^{-\infty}).$$
If $x$ and $y$ belong to different fibres, the phase
$t\psi(x,r_{\eta}y)+\eta $ has no critical point, so
$S_N(x,y)=O(N^{-\infty})$; this estimation is uniform outside a
neighbourhood of $\pi^{-1}(\diag(M))$.

Using the local expression of the phase, one can derive as in
\cite{Shiffman2002} an
expression for $S_N$ at a neighbourhood of size $N^{-1/2}$ of the
diagonal. Let $\Omega_N \subset \C^n\times \R$ be the set of those
$(z,\theta)$ such that $(z/\sqrt{N},\theta/N)$ belongs to the domain
of $\rho$. 
\begin{prop}[\cite{Shiffman2002}, Theorem 3.1]
There exists a sequence $(b_k)_{k\in \N}$ of polynomials on $\R^{4n}$,
such that each $b_k$ is of same parity as $k$, and a smooth function
$R_K$ on $\C^{2n}\times \N$, bounded on the compact sets of $\C^{2n}$
independently of the second variable, such that for all $N$, for
all $(z,w,\theta,phi)\in \Omega_N^2\times \R^2$, there holds
\begin{multline}
  \label{eq:SNexpscaledschiffman}
  N^{-n}e^{i(\phi-\theta)}S_N\left(\rho\left(\frac{z}{\sqrt{N}},\frac{\theta}{N}\right),\left(
    \frac{w}{\sqrt{N}},\frac{
      \phi}{N}\right)\right)\\=\Pi_1(z,w)\left(1+\sum_{k=1}^KN^{-k/2}b_k(z,w,P_0)+N^{-(K+1)/2}R_K(z,w,N)\right)\\+O(N^{-\infty}).
\end{multline}
Here, $\Pi_1$ is the kernel of the projector on the Bargmann space, as in equation (\ref{eq:projbarg}).
\end{prop}
\begin{rem}
  The next step is Proposition \ref{prop:rem}, an estimate for $R_K$ that is valid in all of
  $\Omega_N^2$. For this, we have to keep the
  $O(N^{-\infty})$ term outside. 

In \cite{Shiffman2002}, the
  $O(N^{-\infty})$ term is absorbed into $R_K$, without altering the
  property that $R_K$ is bounded on compact sets independently on
  $N$. However, if an estimate such that the one in Proposition
  \ref{prop:rem} did hold without the supplementary $O(N^{-\infty})$
  term, then one could deduce exponential estimates for the off-diagonal of $S_N$, that is, $|S_N(x,y)|\leq e^{-cN|x-y|^2}$ for some $C$. Such results are indeed known \cite{Berman2008} but cannot be obtained via the Boutet-Sj\"ostrand parametrix because the Boutet-Guillemin construction \cite{BoutetdeMonvel1981} adapts the Szeg\H{o} kernel parametrix to the more general case of almost K\"ahler manifolds, where exponential estimates for the off-diagonal of $S_N$ fail to hold \cite{christ2013upper}.
\end{rem}

The method of proof for the last proposition can be in fact adapted to compute $S_N$ in a fixed
neighbourhood of a point on the diagonal, giving a result close to the Theorem 4.18 of \cite{ma2007holomorphic}, which also appears in \cite{charles2000, Berman2008}. Recall
$$S_N(x,y)=N\iint
\exp(iN(t\psi(x,r_{\eta}y)+\eta))s(x,r_{\eta}y,Nt)\dd t \dd \eta + O(N^{-\infty}).$$
Replacing $\psi$ and $s$ by their expressions we get, after a change
of variables,
\begin{multline*}S_N(\rho(z,\theta),\rho(w,\phi))\\=Ne^{iN(\theta-\phi)}\iint
e^{-N(t(1-A(z,w)e^{i\eta}) -i\eta)}e^{in\eta}\xi(z,w,Nt)\dd t \dd \eta+O(N^{-\infty}).
\end{multline*}
We cannot use the stationary phase lemma, except if $z=w$, because the phase has no critical points. But $\psi $
and $s$ depend holomorphically on $e^{i\eta}$. Thus, we can
replace this integral, which is a contour integral on the unit circle,
with an integral on the circle of radius $|A(z,w)|$ in order to get a phase with a critical point. This
corresponds to formally changing $\eta $ into $\eta -i\log(|A(z,w)|)$
in the computations.
The integral now reads
\begin{multline*}S_N(\rho(z,\theta),\rho(w,\phi))=\\NA(z,w)^Ne^{iN(\theta-\phi)}\iint
e^{-N(t(1-e^{i\eta})-i\eta)}e^{in\eta}\cfrac{\xi(z,w,Nt)}{A(z,w)^n}\dd t \dd \eta+O(N^{-\infty}).\end{multline*}
The last part of the product can now be computed using a stationary
phase lemma, and the fact that $\xi$ is a classical symbol. Hence, we recover a result similar to \cite{ma2007holomorphic,charles2003berezin,Berman2008}:
\begin{prop}
	There exists a neighbourhood $V$ of $(\pi,\pi)^{-1}\diag(M)$ in $X\times X$ such that one has, in local Heisenberg coordinates around a point $P_0\in \diag(X)$ with values in $V$, and for each integer $K$:
\begin{multline}
\label{eq:SNexpfixedSchiffman}
S_N(\rho(z,\theta),\rho(w,\phi))\\=N^{n}e^{iN(\theta-\phi)}A(z,w)^N\left(\sum_{j=0}^KN^{-j}B_j(z,w,P_0)+N^{-(K+1)}r_K(z,w,N,P_0)\right)\\+O(N^{-\infty}).\end{multline}
Each $B_j$ is smooth and $B_0$ is $\frac{1}{\pi^n}$ on the diagonal. Moreover, $r_K$ is bounded in a compact subset of the
domain of definition of $\rho $, independently of $P_0$ and $N$.
\end{prop}
On the diagonal set, $B_0(z,z,P_0)=\cfrac{1}{\pi^n}$ because $S_N$ is
a projector.

Since, in a neighbourhood small enough of the diagonal, one has $$|A(z,w)|\leq 1-\frac 14 |z-w|^2,$$
equation (\ref{eq:SNexpscaledschiffman}) can be deduced from equation (\ref{eq:SNexpfixedSchiffman}). This way, we obtain an estimate on the remainder:

\begin{prop}\label{prop:rem}
	In the equation (\ref{eq:SNexpscaledschiffman}), there exist
        $C$ and $m$ such that the remainder $R_K$ satisfies, for every
        $N$, for every $z$ and $w$ in $\Omega_N$, the inequality:
	\begin{equation*}
		|R_K(z,w,N,P_0)|\leq Ce^{\frac 14|z-w|^2}(1+|z|^m+|w|^m).
	\end{equation*}
\end{prop}
\begin{preuve}
	Rescaling the formula (\ref{eq:SNexpfixedSchiffman}) yields:
	\begin{multline*}
		 N^{-n}e^{i(\phi-\theta)}S_N\left(\rho\left(\frac{z}{\sqrt{N}},\frac{\theta}{N}\right),\left(
		 \frac{w}{\sqrt{N}},\frac{
		 	\phi}{N}\right)\right)\\=A\left(\frac{z}{\sqrt{N}},\frac{w}{\sqrt{N}}\right)^N\left(\sum_{j=0}^KN^{-j}B_j\left(\frac{z}{\sqrt{N}},\frac{w}{\sqrt{N}}\right)+N^{-(K+1)}r_K\left(\frac{z}{\sqrt{N}},\frac{w}{\sqrt{N}},N\right)\right)\\+O(N^{-\infty})
	\end{multline*}
	
	The functions $B_j$ are smooth, and $r_K$ is bounded independently of $N$. Thus, applying a Taylor expansion at the origin, there exist polynomials $b^s_j$, and a function $r^s_K$ with polynomial growth independent of $N$, such that
	\begin{multline}\label{eq:HalfTaylor}
	N^{-n}e^{i(\phi-\theta)}S_N\left(\rho\left(\frac{z}{\sqrt{N}},\frac{\theta}{N}\right),\left(
	\frac{w}{\sqrt{N}},\frac{
		\phi}{N}\right)\right)\\=A\left(\frac{z}{\sqrt{N}},\frac{w}{\sqrt{N}}\right)^N\left(\sum_{j=0}^{2K+1}N^{-j/2}b^s_j(z,w)+N^{-(K+1)}r^s_K(z,w,N)\right)\\+O(N^{-\infty}).
	\end{multline}
	
	Let again $R_A$ be such that $A(z,w)=\pi^n\Pi_1(z,w)e^{R_A(z,w)}$. We wish to control, for any integer $N$, the Taylor expansion at zero of $$g_N:(z,w)\mapsto e^{NR_A\left(\frac{z}{\sqrt{N}},\frac{w}{\sqrt{N}}\right)}.$$
	
	For every multi-index $\alpha$, the derivative of degree $\alpha$ of $g_N$ is a sum of terms of the form $$e^{NR_A\left(\frac{z}{\sqrt{N}},\frac{w}{\sqrt{N}}\right)}\prod_{i=1}^{4n}N^{1-\frac 12|\beta_i|}\partial^{\beta_i}_iR_A\left(\frac{z}{\sqrt{N}},\frac{w}{\sqrt{N}}\right),$$ where each index $\beta_i$ is nonzero and $\sum \beta_i=\alpha$.
	
	Recall that $A$ and $\pi^n\Pi_1$ coincide up to order $2$ at the origin. In particular, the derivatives of order less than $2$ of $R_A$ vanish at the origin. It follows that a term of the form above is nonzero at the origin only if, for each $1\leq i \leq 4n$, there holds $\beta_i\geq 3$. In particular, for each $\alpha$ there holds $$\partial^{\alpha}g_N(0,0)=O(N^{-|\alpha|/6}).$$ Moreover, $\partial^{\alpha}g_N(0,0)$ is always a polynomial in $N^{-1/2}$. 
	
	As we want to write an expansion with a remainder in $O(N^{-K-1})$, let us consider the Taylor expansion of $g_N$ at order $6K+5$. To control the remainder, we make use again of the fact that $R_A$ is smooth on a compact set and that $R_A(z,w)=O(|z|^3,|w|^3)$ at the origin. If $\beta_i=1$, then there is a constant $C$ such that, for every $(z,w)$ and every $N$, one has $$\left|\partial_i^{\beta_i}R_A\left(\frac{z}{\sqrt{N}},\frac{w}{\sqrt{N}}\right) \right|\leq CN^{-1}(|z|^2+|w|^2).$$ Similarly, if $\beta_i=2$, there exists a constant $C$ such that, for every $(z,w)$ and every $N$, one has
	$$\left|\partial_i^{\beta_i}R_A\left(\frac{z}{\sqrt{N}},\frac{w}{\sqrt{N}}\right)\right|
        \leq CN^{-1/2}(|z|+|w|).$$ If $\beta_i\geq 3$ we simply use
        the fact that the function $\partial_i^{\beta_i}R_A$ is bounded
        on its set of definition. It follows that for every $\alpha$
        there exist $m$ and $C$ such that, for every $N$, for every
        $z,w \in \Omega_N$, one has
	$$\left|\partial^{\alpha}g_N(z,w)\right| \leq CN^{-|\alpha|/6}(1+|z|^m+|w|^m)\left|g_N(z,w)\right|. $$
	
	Recall now from Proposition \ref{prop:remainder} that $$|g_1(z,w)| \leq e^{\frac 14 |z-w|^2}.$$ From the definition of $g_N$ one deduces that $$|g_N(z,w)| \leq e^{\frac 14 |z-w|^2}.$$ Thus the Taylor expansion of $g_N$ of order $6K+5$ at the origin takes the following form:
	$$g_N(z,w)=\sum_{j=0}^{2K+1}N^{-j/2}b^{\psi}_j(z,w) + N^{-K-1}r^{\psi}_K(z,w,N).$$ 
	Here, the $b^{\psi}_j$ are polynomials, and there exist $C$ and $m$ such that, for every $z,w$ and every $N$, one has
	$$|r^{\psi}_K(z,w,N)|\leq (1+|z|^m+|w|^m)e^{-\frac 14 |z-w|^2}.$$
	
	We now return to equation (\ref{eq:HalfTaylor}). Replacing $A$ with $\pi^n\Pi_1e^{R_A}$, using the previous expression of $g_N$ and expanding, we find equation (\ref{eq:SNexpscaledschiffman}) with the desired control of $R_K$. 
\end{preuve}

\bibliographystyle{abbrv}
\bibliography{math}
\end{document}